\documentclass[12pt]{article}
\usepackage{amssymb}

\newtheorem{definition}{Definition}[section]
\newtheorem{theorem}[definition]{Theorem}
\newtheorem{lemma}[definition]{Lemma}
\newtheorem{corollary}[definition]{Corollary}
\newtheorem{remark}[definition]{Remark}

\newtheorem{note}[definition]{Note}

\def\K{\mathbb K}

\def\fld{\mathbb K}
\def\Mdf{{\mbox{Mat}}_{d+1}(\K)}
\def\alg{\mathcal A}

\begin{document}

\title{ \bf Leonard pairs and the $q$-Racah polynomials}
\author{Paul Terwilliger  
}
\date{}
\maketitle
\begin{abstract} 
Let $\K$ denote a field, and let $V$ denote a  
vector space over $\K$ with finite positive dimension.
We consider a pair
of linear transformations
$A:V\rightarrow V$ and $A^*:V\rightarrow V$
that satisfy the following two conditions:
\begin{enumerate}
\item There exists a basis for $V$ with respect to which
the matrix representing $A$ is irreducible tridiagonal and the matrix
representing $A^*$ is diagonal.
\item There exists a basis for $V$ with respect to which
the matrix representing $A$ is diagonal and the matrix
representing $A^*$ is irreducible tridiagonal.
\end{enumerate}

\medskip
\noindent
We call such a pair a {\it Leonard pair} on $V$. 
In the appendix to 
\cite{LS99} we outlined a correspondence between
Leonard pairs and a class of orthogonal polynomials
consisting of the $q$-Racah polynomials and some related
polynomials of the Askey scheme. 
We also outlined how, for the polynomials in this class,
the
3-term recurrence,  difference equation, Askey-Wilson
duality, and orthogonality can be obtained in a uniform 
manner from 
the corresponding Leonard pair.
The purpose of this paper is to
 provide proofs for the assertions which we made in 
that appendix. 
\end{abstract}

\section{Leonard pairs}
We begin by recalling the notion of a {\it Leonard pair} 
\cite{TD00},
\cite{LS99},
\cite{qSerre},
\cite{LS24},
\cite{conform},
\cite{lsint},
\cite{Terint},
\cite{TLT:split},
\cite{TLT:array}.
We will use the following terms. Let $X$ denote a square matrix.
Then $X$ is called {\it tridiagonal} whenever each nonzero entry
lies on either the diagonal, the subdiagonal, or the superdiagonal.
Assume $X$ is tridiagonal. Then $X$ is called {\it irreducible}\
whenever each entry on the subdiagonal is nonzero and each entry
on the superdiagonal is nonzero.

\medskip
\noindent We now define a Leonard pair. For the rest of this paper
$\K$ will denote a field.

\begin{definition} \cite{LS99}
\label{def:lprecall}
\rm
Let 
 $V$ denote a  
vector space over $\K$ with finite positive dimension.
By a {\it Leonard pair} on $V$,
we mean an ordered pair of linear transformations
$A:V\rightarrow V$ and $A^*:V\rightarrow V$ 
that
 satisfy both (i), (ii) below. 
\begin{enumerate}
\item There exists a basis for $V$ with respect to which
the matrix representing $A$ is irreducible tridiagonal and the matrix
representing $A^*$ is diagonal.
\item There exists a basis for $V$ with respect to which
the matrix representing $A$ is diagonal and the matrix
representing $A^*$ is irreducible tridiagonal.
\end{enumerate}
\end{definition}

\begin{note}
\rm
According to a common notational convention
$A^*$ denotes the conjugate-transpose of $A$.
We are not using this convention. In a Leonard
pair $A,A^*$, the linear transformations $A$ and $A^*$
are arbitrary subject
to (i),  (ii) above.
\end{note}

%
%
%
%
\medskip

\section{An example}

Here is an example of a Leonard pair.
Set 
$V={\K}^4$ (column vectors), set 
\begin{eqnarray*}
A = 
\left(
\begin{array}{ c c c c }
0 & 3  &  0    & 0  \\
1 & 0  &  2   &  0    \\
0  & 2  & 0   & 1 \\
0  & 0  & 3  & 0 \\
\end{array}
\right), \qquad  
A^* = 
\left(
\begin{array}{ c c c c }
3 & 0  &  0    & 0  \\
0 & 1  &  0   &  0    \\
0  & 0  & -1   & 0 \\
0  & 0  & 0  & -3 \\
\end{array}
\right),
\end{eqnarray*}
and view $A$ and $A^*$  as linear transformations from $V$ to $V$.
We assume 
the characteristic of $\K$ is not 2 or 3, to ensure
$A$ is irreducible.
Then $A, A^*$ is a Leonard
pair on $V$. 
Indeed, 
condition (i) in Definition
\ref{def:lprecall}
is satisfied by the basis for $V$
consisting of the columns of the 4 by 4 identity matrix.
To verify condition (ii), we display an invertible  matrix  
$P$ such that 
$P^{-1}AP$ is 
diagonal and 
$P^{-1}A^*P$ is
irreducible tridiagonal.
Set 
\begin{eqnarray*}
P = 
\left(
\begin{array}{ c c c c}
1 & 3  &  3    &  1 \\
1 & 1  &  -1    &  -1\\
1  & -1  & -1  & 1  \\
1  & -3  & 3  & -1 \\
\end{array}
\right).
\end{eqnarray*}
 By matrix multiplication $P^2=8I$, where $I$ denotes the identity,   
so $P^{-1}$ exists. Also by matrix multiplication,    
\begin{equation}
AP = PA^*.
\label{eq:apeq}
\end{equation}
Apparently
$P^{-1}AP$ is equal to $A^*$ and is therefore diagonal.
By (\ref{eq:apeq}) and since $P^{-1}$ is
a scalar multiple of $P$, we find
$P^{-1}A^*P$ is equal to $A$ and is therefore irreducible tridiagonal.  Now 
condition (ii) of  Definition 
\ref{def:lprecall}
is satisfied
by the basis for $V$ consisting of the columns of $P$. 

\medskip
\noindent The above example is a member of the following infinite
family of Leonard pairs.
For any nonnegative integer $d$  
the pair
\begin{equation}
A = 
\left(
\begin{array}{ c c c c c c}
0 & d  &      &      &   &{\bf 0} \\
1 & 0  &  d-1   &      &   &  \\
  & 2  &  \cdot    & \cdot  &   & \\
  &   & \cdot     & \cdot  & \cdot   & \\
  &   &           &  \cdot & \cdot & 1 \\
{\bf 0} &   &   &   & d & 0  
\end{array}
\right),
\qquad A^*= \hbox{diag}(d, d-2, d-4, \ldots, -d)
\label{eq:fam1}
\end{equation}
is a Leonard pair on the vector space $\K^{d+1} $,
provided the 
 characteristic of $\K$ is zero or an odd prime greater than $d$.
This can be  proved by modifying the 
 proof for $d=3$ given above. One shows  
$P^2=2^dI$  and $AP= PA^*$, where 
$P$ denotes the matrix with $ij$ entry
\begin{equation}
P_{ij} =  
\Biggl({{ d }\atop {j}}\Biggr) {{}_2}F_1\Biggl({{-i, -j}\atop {-d}}
\;\Bigg\vert \;2\Biggr)
\qquad \qquad (0 \leq i,j\leq d)
\label{eq:ex1}
\end{equation}
\cite[Section 16]{LS24}.
We follow the standard notation for
hypergeometric series \cite[p. 3]{GR}. 

\medskip

\section{Leonard systems}
\medskip
\noindent When working with a Leonard pair, 
it is often convenient to consider a closely related
and somewhat more abstract object called
a {\it Leonard system}.
In order to define this we first make an observation about
Leonard pairs.

\begin{lemma} \cite[Lemma 1.3]{LS99}
\label{lem:preeverythingtalkS99}
Let $V$ denote a vector space over $\K$ with finite positive
dimension and 
let $A, A^*$ denote a Leonard pair on $V$. Then
the eigenvalues of $A$ are mutually distinct
and contained in $\K$.
Moreover, the eigenvalues of
$A^*$ are mutually distinct and contained in $\fld$.
\end{lemma}

\noindent To prepare for our definition of  a Leonard system,
we recall a few concepts from  linear algebra. 
Let $d$  denote  a  nonnegative
integer 
and let 
$\Mdf$ denote the $\fld$-algebra consisting of all
$d+1$ by $d+1$ matrices that have entries in $\K$. We
index the rows and columns by $0,1,\ldots, d$.
We let ${\K}^{d+1}$ denote the $\K$-vector space consisting
of all $d+1$ by $1$ matrices that have entries in $\K$.
We index the rows by $0,1,\ldots, d$. We view 
 ${\K}^{d+1}$ as a left module for 
$\Mdf$.
We observe this module is irreducible.
For the rest of this paper we  
let  
$\mathcal A$
denote a $\fld$-algebra
isomorphic to 
$\hbox{Mat}_{d+1}(\K)$. 
When we refer to an $\mathcal A$-module we mean a left 
 $\mathcal A$-module.
Let $V$ denote an irreducible
$\mathcal A$-module.
We remark that $V$ is unique up to isomorphism of 
$\mathcal A$-modules, and 
that $V$ has dimension $d+1$.
Let $v_0, v_1, \ldots, v_d$ denote a basis for $V$.
For $X \in \mathcal A$ and
$Y \in 
\mbox{Mat}_{d+1}(\K)$,
we say 
$Y$ {\it represents $X$ with respect to $v_0, v_1, \ldots, v_d$}
whenever
$Xv_j = \sum_{i=0}^d Y_{ij}v_i$ for 
$0 \leq j\leq d$.
Let $A$ denote an element of $\mathcal A$.
We say $A$ is 
{\it multiplicity-free} whenever it has $d+1$ 
mutually distinct  eigenvalues in $\K$.
Let $A$ denote a  multiplicity-free element of $\alg$.
Let $\theta_0, \theta_1, \ldots, \theta_d$ denote an ordering of 
the eigenvalues
of $A$, and for $0 \leq i \leq d$   put 
\begin{equation}
E_i = \prod_{{0 \leq  j \leq d}\atop
{j\not=i}} {{A-\theta_j I}\over {\theta_i-\theta_j}},
\label{eq:primiddef}
\end{equation}
where $I$ denotes the identity of $\cal A$.
We observe
(i) 
$AE_i = \theta_iE_i \; (0 \leq i \leq d)$;
(ii) $E_iE_j = \delta_{ij}E_i \;(0 \leq i,j\leq d)$;
(iii) $\sum_{i=0}^d E_i = I $;
(iv) $A = \sum_{i=0}^d \theta_i E_i.$
Let $\mathcal D$ denote the subalgebra of $\mathcal A$ generated
by $A$. Using (i)--(iv) we find the sequence
 $E_0, E_1, \ldots, E_d$ is a  basis for the
$\K$-vector space $\cal D$.
We call $E_i$  the {\it primitive idempotent} of
$A$ associated with $\theta_i$.
It is helpful to think of these primitive idempotents as follows. 
Observe 
\begin{eqnarray*}
V = E_0V + E_1V + \cdots + E_dV \qquad \qquad (\mbox{direct sum}).
\end{eqnarray*}
For $0\leq i \leq d$, $E_iV$ is the (one dimensional) eigenspace of
$A$ in $V$ associated with the 
eigenvalue $\theta_i$, 
and $E_i$ acts  on $V$ as the projection onto this eigenspace. 
We remark that $\lbrace A^i |0 \leq i \leq d\rbrace $ 
is a basis for the 
$\K$-vector space $\mathcal D$ and that
$\prod_{i=0}^d (A-\theta_iI)=0$. By a {\it Leonard pair in ${\mathcal A}$}
we mean an ordered pair of elements taken from $\mathcal A$ that act
on $V$ as a Leonard pair in the sense of Definition
\ref{def:lprecall}. We call $\mathcal A$ the {\it ambient algebra} of
the pair and say the pair is {\it over $\K$}. 
We refer to $d$ as the {\it diameter} of the pair.
We now define a Leonard system.

\begin{definition} \cite{LS99}
\label{def:defls}
\label{eq:ourstartingpt}
\rm
By a {\it Leonard system} in $\mathcal A$ we mean a 
sequence
$\Phi:=(A;A^*; \lbrace E_i\rbrace_{i=0}^d; $ $ 
\lbrace E^*_i\rbrace_{i=0}^d)$
that satisfies (i)--(v) below. 
\begin{enumerate}
\item Each of $A,A^*$ is a multiplicity-free element in $\mathcal A$.
\item $E_0,E_1,\ldots,E_d$ is an ordering of the primitive 
idempotents of $A$.
\item $E^*_0,E^*_1,\ldots,E^*_d$ is an ordering of the primitive 
idempotents of $A^*$.
\item ${\displaystyle{
E_iA^*E_j = \cases{0, &if $\;\vert i-j\vert > 1$;\cr
\not=0, &if $\;\vert i-j \vert = 1$\cr}
\qquad \qquad 
(0 \leq i,j\leq d)}}$.
\item ${\displaystyle{
 E^*_iAE^*_j = \cases{0, &if $\;\vert i-j\vert > 1$;\cr
\not=0, &if $\;\vert i-j \vert = 1$\cr}
\qquad \qquad 
(0 \leq i,j\leq d).}}$
\end{enumerate}
We refer to $d$ as the {\it diameter} of $\Phi$ and say 
$\Phi$ is {\it over } $\K$.  We call $\mathcal A$
the {\it ambient algebra} of $\Phi$. 
\end{definition}

\noindent We comment on how 
Leonard pairs and Leonard systems are related.
In the following discussion $V$ denotes
an irreducible $\mathcal A$-module. Let 
$(A;A^*; \lbrace E_i\rbrace_{i=0}^d;  
\lbrace E^*_i\rbrace_{i=0}^d) $
denote a Leonard system in $\mathcal A$. 
For $0 \leq i \leq d$ let $v_i$ denote a nonzero vector in 
$E_iV$. Then the sequence 
$v_0, v_1, \ldots, v_d$ is a basis for $V$ that satisfies
Definition 
\ref{def:lprecall}(ii).
For $0 \leq i \leq d$ let $v^*_i$ denote a nonzero vector in 
$E^*_iV$. Then the sequence 
$v^*_0, v^*_1, \ldots, v^*_d$ is a basis for $V$ which satisfies
Definition 
\ref{def:lprecall}(i). By these comments the pair $A,A^*$ is
a Leonard pair in $\mathcal A$. Conversely let $A,A^*$ denote
a Leonard pair in $\mathcal A$. Then each of $A, A^*$ is multiplicity-free
by Lemma
\ref{lem:preeverythingtalkS99}.
 Let 
$v_0, v_1, \ldots, v_d$ denote a basis for $V$ that satisfies
Definition 
\ref{def:lprecall}(ii). For $0 \leq i \leq d$ the vector
$v_i$ is an eigenvector for $A$;
 let
$E_i$ denote the corresponding primitive idempotent.
 Let 
$v^*_0, v^*_1, \ldots, v^*_d$ denote a basis for $V$ that satisfies
Definition 
\ref{def:lprecall}(i). For $0 \leq i \leq d$ the vector
$v^*_i$ is an eigenvector for $A^*$;
 let
$E^*_i$ denote the corresponding primitive idempotent.
Then 
$(A;A^*; \lbrace E_i\rbrace_{i=0}^d;  
\lbrace E^*_i\rbrace_{i=0}^d)$
is a Leonard system in $\mathcal A$.
In summary we have the following.

\begin{lemma} Let $A$ and $A^*$ denote elements of $\mathcal A$.
Then the pair $A,A^*$ is a Leonard pair in $\mathcal A$ if and only 
if the following (i), (ii) hold.
\begin{enumerate}
\item Each of $A,A^*$ is multiplicity-free.
\item There exists an ordering $E_0, E_1, \ldots, E_d$ of the primitive
idempotents of $A$ and there exists
an ordering $E^*_0, E^*_1, \ldots, E^*_d$ of the primitive
idempotents of $A^*$ such that
$(A;A^*; $ $\lbrace E_i\rbrace_{i=0}^d;  
\lbrace E^*_i\rbrace_{i=0}^d)$ is a Leonard system in $\mathcal A$.
\end{enumerate}
\end{lemma}

\medskip
\noindent We recall the notion of {\it isomorphism} for Leonard pairs and
Leonard systems.

\begin{definition} 
\rm
Let $A,A^*$ and $B,B^*$ denote Leonard pairs over
$\K$. By an {\it isomorphism of Leonard pairs from $A,A^*$ to $B,B^*$} 
we mean an isomorphism of $\K$-algebras from the ambient algebra of
$A,A^*$ to the ambient algebra of $B,B^*$ that sends $A$ to $B$
and $A^*$ to $B^*$.  The Leonard pairs $A,A^*$ and $B,B^*$ are
said to be {\it isomorphic} whenever there exists an isomorphism
of Leonard pairs from 
 $A,A^*$ to  $B,B^*$.
\end{definition} 
\noindent 
Let $\Phi$ 
denote the Leonard system from Definition
\ref{eq:ourstartingpt}
and let 
$\sigma :\alg \rightarrow {\cal A}'$ denote an isomorphism of
$\K$-algebras. We write 
$\Phi^{\sigma}:= 
(A^\sigma;A^{*\sigma};\lbrace E^\sigma_i\rbrace_{i=0}^d;  
\lbrace E^{*\sigma}_i\rbrace_{i=0}^d)$
and observe 
$\Phi^{\sigma}$ 
is a Leonard  system in ${\cal A }'$.

\begin{definition}
\label{def:isolsS99o}
\rm
Let $\Phi$ 
and  
 $\Phi'$ 
denote Leonard systems over $\K$.
 By an {\it isomorphism of Leonard  systems
 from $\Phi $ to $\Phi'$} we mean an isomorphism 
of $\K $-algebras
$\sigma $
from the ambient algebra of $\Phi$ to the ambient algebra of
$\Phi'$ such 
that  $\Phi^\sigma = \Phi'$. 
The Leonard systems $\Phi $, $\Phi'$
are said to be {\it isomorphic} whenever there exists
an isomorphism of Leonard  systems from $\Phi $ to $\Phi'$. 
\end{definition}

\noindent We have a remark.
Let  $\sigma :\mathcal A \rightarrow 
\mathcal A$ denote  any map.
By the Skolem-Noether theorem 
\cite[Corollary 9.122]{rotman},
$\sigma $ is an isomorphism of $\K$-algebras
if and only if there exists an invertible $S \in {\mathcal A}$ such that
$X^\sigma = S X S^{-1}$ for all  $X \in  {\mathcal A}$.

%

\section{The $D_4$ action}

\medskip
\noindent A given Leonard system  can be modified in  several
ways to get a new Leonard system. For instance, 
let $\Phi$ 
 denote the Leonard system from 
Definition 
\ref{eq:ourstartingpt}.
Then each of the following three sequences is a Leonard system
in $\mathcal A$.
\begin{eqnarray*}
 \;\Phi^*&:=& (A^*; A; \lbrace E^*_i\rbrace_{i=0}^d; 
 \lbrace E_i\rbrace_{i=0}^d),
\\
\Phi^{\downarrow}&:=& (A; A^*; \lbrace E_i\rbrace_{i=0}^d ;
\lbrace E^*_{d-i} \rbrace_{i=0}^d),
\\
\Phi^{\Downarrow} 
&:=& (A;A^*; \lbrace E_{d-i}\rbrace_{i=0}^d; \lbrace E^*_i\rbrace_{i=0}^d).
\end{eqnarray*}
Viewing $*, \downarrow, \Downarrow$
as permutations on the set of all Leonard systems,
\begin{eqnarray}
&&\qquad \qquad \qquad  *^2 \;=\;  
\downarrow^2\;= \;
\Downarrow^2 \;=\;1,
\qquad \quad 
\label{eq:deightrelationsAS99}
\\
&&\Downarrow *\; 
=\;
* \downarrow,\qquad \qquad   
\downarrow *\; 
=\;
* \Downarrow,\qquad \qquad   
\downarrow \Downarrow \; = \;
\Downarrow \downarrow.
\qquad \quad 
\label{eq:deightrelationsBS99}
\end{eqnarray}
The group generated by symbols 
$*, \downarrow, \Downarrow $ subject to the relations
(\ref{eq:deightrelationsAS99}),
(\ref{eq:deightrelationsBS99})
is the dihedral group $D_4$.  
We recall $D_4$ is the group of symmetries of a square,
and has 8 elements.
Apparently $*, \downarrow, \Downarrow $ induce an action of 
 $D_4$ on the set of all Leonard systems.
%
%

\medskip
\noindent For the rest of this paper we will use the following
notational convention.

\begin{definition}
\label{def:notconv}
\rm
Let $\Phi$ denote a Leonard system. For any element $g$ in
the group $D_4$ and for any object $f$ that we associate with
$\Phi$, we let $f^g$ denote the corresponding object for the
Leonard system 
$\Phi^{g^{-1}}$. We have been using this
convention all along; an example is $E^*_i(\Phi)= 
E_i(\Phi^*)$.
\end{definition}

\section{The structure of a Leonard system}

\noindent In this section we establish a few basic 
facts concerning Leonard systems.
We begin with a definition and two routine lemmas.

\begin{definition}
\label{def:evseq}
\rm
Let $\Phi$ denote the Leonard system  from  
Definition \ref{eq:ourstartingpt}.
For $0 \leq i \leq d$, 
we let $\theta_i $ (resp. $\theta^*_i$) denote the eigenvalue
of $A$ (resp. $A^*$) associated with $E_i$ (resp. $E^*_i$).
We refer to  $\theta_0, \theta_1, \ldots, \theta_d$ as the 
{\it eigenvalue sequence} of $\Phi$.
We refer to  $\theta^*_0, \theta^*_1, \ldots, \theta^*_d$ as the 
{\it dual eigenvalue sequence} of $\Phi$. We observe 
 $\theta_0, \theta_1, \ldots, \theta_d$ are mutually distinct
 and contained in $\K$. Similarly
  $\theta^*_0, \theta^*_1, \ldots, \theta^*_d$  
 are mutually distinct
 and contained in $\K$. 
\end{definition}

\begin{lemma}
\label{lem:step0}
Let $\Phi$ denote the Leonard system
from Definition  
\ref{eq:ourstartingpt} and let $V$ denote an irreducible
$\mathcal A$-module. For $0 \leq i \leq d$ let $v_i$
denote a nonzero vector in $E^*_iV$ and observe
$v_0, v_1, \ldots, v_d$ is a basis for $V$.
Then  (i), (ii) hold below.
\begin{enumerate}
\item
For $0 \leq i \leq d$ the
matrix in 
$\mbox{Mat}_{d+1}(\K)$
that represents $E^*_i$ with respect to
$v_0, v_1, \ldots, v_d$ has $ii$ entry 1 and all other entries 0.
\item
The matrix in 
$\mbox{Mat}_{d+1}(\K)$
that represents $A^*$ with respect to
$v_0, v_1, \ldots, v_d$ is equal to
$\mbox{diag}(\theta^*_0, \theta^*_1, \ldots, \theta^*_d)$.
\end{enumerate}
\end{lemma}

\begin{lemma}
\label{lem:step1}
Let $A$ denote an irreducible tridiagonal matrix
in 
$\mbox{Mat}_{d+1}(\K)$.
Pick any integers $i,j$ $(0 \leq i,j\leq d)$.
Then (i)--(iii) hold below.
\begin{enumerate}
\item
The entry $(A^r)_{ij}= 0$ if  $r<\vert i-j \vert,
\qquad \qquad (0 \leq r \leq d).$
\item
Suppose $i\leq j$. Then 
the entry
$(A^{j-i})_{ij} = \prod_{h=i}^{j-1}A_{h,h+1}.$
Moreover 
$(A^{j-i})_{ij} \not=0$.
\item
Suppose $i\geq j$. Then 
the entry
$(A^{i-j})_{ij} = \prod_{h=j}^{i-1}A_{h+1,h}.$
Moreover 
$(A^{i-j})_{ij} \not=0$.
\end{enumerate}
\end{lemma}

\begin{theorem}
\label{eq:lsmatbasis}
Let $\Phi$ denote the Leonard system
from Definition  
\ref{eq:ourstartingpt}. Then the elements
\begin{equation}
A^rE^*_0A^s \qquad \qquad (0 \leq r,s\leq d)
\label{eq:lpbasis}
\end{equation}
form a basis for the $\K$-vector space $\cal A$.
\end{theorem}

\noindent {\it Proof:}
The number of elements in 
(\ref{eq:lpbasis}) is equal to $(d+1)^2$, and this number is the dimension
of $\cal A$. Therefore
it suffices
to show 
the elements
in (\ref{eq:lpbasis})
are linearly independent. To do this,
we represent 
the elements 
in (\ref{eq:lpbasis}) by matrices. 
Let $V$ denote an irreducible 
$\mathcal A$-module.
For $0 \leq i \leq d$ let $v_i$ denote a nonzero vector in 
$E^*_iV$, and observe $v_0, v_1,\ldots, v_d$ is a basis for
$V$. For the purpose of this proof, let us identify each element
of $\mathcal A$ with the matrix in 
$\hbox{Mat}_{d+1}(\K)$ that 
represents it with respect to the
basis $v_0, v_1, \ldots, v_d$. Adopting this point of view
we find 
$A$ is irreducible tridiagonal and $A^*$ is diagonal. 
For $0 \leq r,s\leq d$ we show  the
entries  of $A^rE^*_0A^s$ satisfy 
\begin{equation}
(A^rE^*_0A^s)_{ij} = 
 \cases{0, 
&$\qquad $if $\quad i>r \quad$ or $\quad j>s $ ;\cr
\not=0,  & $\qquad $if $\quad i=r\quad $ and $\quad j=s$\cr} 
\qquad \qquad (0 \leq i,j\leq d).
\label{eq:lowert}
\end{equation}
By Lemma
\ref{lem:step0}(i)
the matrix $E^*_0$ has $00$ entry 1 and all other entries $0$.
Therefore 
\begin{equation}
(A^rE^*_0A^s)_{ij} =  (A^r)_{i0} (A^s)_{0j} \qquad \qquad (0 \leq i,j\leq d).
\label{eq:redent}
\end{equation}
We mentioned $A$ is irreducible tridiagonal.
Applying Lemma \ref{lem:step1} we find that for 
$0 \leq i \leq d$ the entry $(A^r)_{i0}$ is zero if $i>r$, and  nonzero if
$i=r$. Similarly for $0 \leq j\leq d$ the entry
 $(A^s)_{0j}$ is zero 
if $j>s$, and nonzero if $j=s$. Combining these facts with 
(\ref{eq:redent}) we routinely obtain
(\ref{eq:lowert}) and it follows  the elements 
(\ref{eq:lpbasis}) are linearly independent.
Apparently the elements 
(\ref{eq:lpbasis}) form a basis  for $\cal A $, as desired.
\hfill $\Box $ \\

\begin{corollary} 
\label{cor:genset}
Let $\Phi$ denote the Leonard system
from Definition 
\ref{eq:ourstartingpt}. Then 
the elements 
$A, E^*_0$ together generate $\cal A$. Moreover
the elements $A,A^*$ together generate $\cal A$.
\end{corollary}

\noindent {\it Proof:}
The first assertion is immediate from Theorem 
\ref{eq:lsmatbasis}.  
The second assertion follows from the first assertion and
the observation that $E^*_0$ is a polynomial in $A^*$. 
\hfill $\Box $ \\

\noindent The following is immediate from Corollary
\ref{cor:genset}.

\begin{corollary} 
\label{cor:gensetlp}
Let $A,A^*$ denote a Leonard pair in 
$\mathcal A$. Then the 
elements $A,A^*$ together generate $\mathcal A$.
\end{corollary}

\noindent We mention a few implications of Theorem 
\ref{eq:lsmatbasis} that will be useful later in the paper.

\begin{lemma}
\label{lem:altbase}
Let $\Phi$ denote the Leonard system
from Definition 
\ref{eq:ourstartingpt}. Let $\cal D$ denote
the subalgebra of $\cal A$ generated by $A$.
 Let $X_0, X_1, \ldots, X_d$
denote a basis for the $\K$-vector space 
$\mathcal D$. Then the elements
\begin{equation}
X_r E^*_0 X_s \qquad \qquad (0 \leq r,s\leq d)
\label{altbase}
\end{equation}
form a basis for the $\K$-vector space 
$\mathcal A$.

\end{lemma}

\noindent {\it Proof:}
 The number of elements in    
(\ref{altbase}) is equal to  $(d+1)^2$, 
and this number is the dimension
of $\cal A$. Therefore it suffices to show the elements
(\ref{altbase}) span $\cal A$. 
But this is immediate from Theorem 
\ref{eq:lsmatbasis}, and since 
each  element 
in (\ref{eq:lpbasis}) is contained in the span of
the elements (\ref{altbase}).
\hfill $\Box $ \\

\begin{corollary}
\label{cor:eibasis}
Let $\Phi$ denote the Leonard system
from Definition
\ref{eq:ourstartingpt}.
Then the elements
\begin{equation}
E_r E^*_0 E_s \qquad \qquad (0 \leq r,s\leq d)
\label{eq:eibase}
\end{equation}
form a basis for the $\K$-vector space $\mathcal A$.
\end{corollary}

\noindent {\it Proof:}
Immediate from  
Lemma \ref{lem:altbase}, with  $X_i = E_i$ for $0 \leq i \leq d$.
\hfill $\Box $ \\

\begin{lemma}
\label{lem:oneiszero}
Let $\Phi$ denote the Leonard system
from Definition 
\ref{eq:ourstartingpt}. Let $\cal D$ denote
the subalgebra of $\cal A$ generated by $A$.
Let $X$ and $Y$ denote elements in $\mathcal D$
and assume $XE^*_0Y=0$. Then $X=0$ or $Y=0$.
\end{lemma}
\noindent {\it Proof:}
Let $X_0, X_1, \ldots, X_d$ denote a basis for 
the $\K$-vector space $\mathcal D$.
Since $X \in {\mathcal D}$ there exists
$\alpha_i \in \K$ $(0 \leq i \leq d)$
such that $X=\sum_{i=0}^d \alpha_i X_i$.
Similarly there exists 
$\beta_i \in \K$
$(0 \leq i \leq d)$
such that $Y=\sum_{i=0}^d \beta_i X_i$.
Evaluating 
$0=XE^*_0Y$ using these equations we get 
 $0= \sum_{i=0}^d \sum_{j=0}^d \alpha_i\beta_j X_iE^*_0X_j$.
From this and 
 Lemma
\ref{lem:altbase}
we find $\alpha_i \beta_j=0$ for $0 \leq i,j\leq d$.
We assume $X\not=0$ and show $Y=0$. 
Since $X \not=0$ there exists an integer $i$
$(0 \leq i \leq d)$ such that $\alpha_i\not=0$.
Now for $0 \leq j\leq d$ we have
$\alpha_i\beta_j=0$ so $\beta_j=0$. It follows
$Y=0$.
\hfill $\Box $ \\

\noindent We finish this section with a comment.

\begin{lemma}
\label{lem:absval}
Let $\Phi$ denote the Leonard system
from Definition 
\ref{eq:ourstartingpt}.
Pick any integers $i,j$ $(0 \leq i,j\leq d)$. Then (i)--(iv) hold
below.
\begin{enumerate}
\item
$
E^*_iA^rE^*_j = 0 \quad \hbox{if } \quad r< \vert i-j \vert,
\qquad (0 \leq r \leq d).
$
\item
Suppose $ i \leq j$. Then
\begin{eqnarray}
\label{eq:leftside}
E^*_i A^{j-i}E^*_j = E^*_iAE^*_{i+1}A \cdots E^*_{j-1}AE^*_j.
\end{eqnarray}
Moreover 
$E^*_i A^{j-i}E^*_j\not=0$.
\item Suppose 
$ i \geq j$. Then
\begin{eqnarray}
\label{eq:rightside}
E^*_i A^{i-j}E^*_j = E^*_iAE^*_{i-1}A \cdots E^*_{j+1}AE^*_j.
\end{eqnarray}
Moroever
$E^*_i A^{i-j}E^*_j\not=0$.
\item Abbreviate $r=|i-j|$.
Then $E^*_iA^rE^*_j$ is a  basis for
the $\K$-vector space 
$E^*_i{\mathcal A}E^*_j$.
\end{enumerate}
\end{lemma}
\noindent {\it Proof:}
Represent the elements of $\Phi$ by matrices 
as in the proof of Theorem
\ref{eq:lsmatbasis}, and use
Lemma
\ref{lem:step1}.
\hfill $\Box $ \\


\section{The antiautomorphism $\dagger$}

We recall the notion of an {\it antiautomorphism} of $\mathcal A$.
Let $\gamma : {\mathcal A} \rightarrow {\mathcal A}$ denote any map.
We call $\gamma $
an {\it antiautomorphism} of $\mathcal A$
whenever $\gamma$ is an
isomorphism of $\K$-vector spaces  
and  $(XY)^\gamma = Y^\gamma X^\gamma $ for
all $X, Y \in {\mathcal A}$.
For example assume ${\mathcal A} = 
\hbox{Mat}_{d+1}(\K)$.
Then $\gamma $ is an antiautomorphism of $\mathcal A$ if and only 
if there exists an invertible element $R$ in ${\mathcal A}$ such that
$X^\gamma = R^{-1}X^t R$ for all $X \in {\mathcal A}$, where
$t$ denotes transpose.
This follows from the Skolem-Noether theorem \cite[Corollary 9.122]{rotman}.

\begin{theorem}
\label{thm:dagger}
Let $A,A^*$ denote a Leonard pair in $\mathcal A$.
Then there exists a unique antiautomorphism $\dagger$ of $\mathcal A$
such that $A^\dagger = A$ and $A^{*\dagger}=A^*$.
Moreover $X^{\dagger \dagger}=X$
for all $X \in {\mathcal A}$.
\end{theorem}
\noindent {\it Proof:}
Concerning existence, 
let $V$ denote an irreducible $\mathcal A$-module. 
By 
Definition
\ref{def:lprecall}(i) there exists a basis
for $V$ with respect to which the matrix representing
$A$ is irreducible tridiagonal and the matrix representing
$A^*$ is diagonal.
Let us denote this basis by  
$v_0, v_1, \ldots, v_d$.
For $X \in \mathcal A$ let $X^\sigma$ denote the matrix
in 
$\mbox{Mat}_{d+1}(\K)$ that represents $X$
with respect to
the basis
$v_0, v_1, \ldots, v_d$.
We observe
$\sigma : {\mathcal A}\rightarrow 
\mbox{Mat}_{d+1}(\K)$ 
is an isomorphism of $\K$-algebras.
We abbreviate $B=A^\sigma$ and observe
$B$ is irreducible tridiagonal.
We abbreviate $B^*=A^{*\sigma}$ and observe
$B^*$ is diagonal.
Let $D$ denote the diagonal matrix in  
$\mbox{Mat}_{d+1}(\K)$ 
that has $ii$ entry
\begin{eqnarray*}
D_{ii} = \frac{B_{01}B_{12} \cdots B_{i-1,i}}
{B_{10}B_{21} \cdots B_{i,i-1}}
\qquad \qquad (0 \leq i \leq d).
\end{eqnarray*}
It is routine to verify $D^{-1}B^tD=B$. Each of $D, B^*$ is
diagonal so $DB^*=B^*D$; also
$B^{*t}=B^*$ so $D^{-1}B^{*t}D=B^*$. Let
$
\gamma :
\mbox{Mat}_{d+1}(\K)
\rightarrow
\mbox{Mat}_{d+1}(\K)
$
denote the map that satisfies
$X^{\gamma} = D^{-1}X^tD$ for all $X \in 
\mbox{Mat}_{d+1}(\K)$.
We observe $\gamma $ is an antiautomorphism
of 
$\mbox{Mat}_{d+1}(\K)$ such that
$B^\gamma=B$ and $B^{*\gamma}=B^*$.
We define the map
$
\dagger :{ \mathcal A}
\rightarrow 
{ \mathcal A}
$
to be the  composition
$\dagger = \sigma \gamma \sigma^{-1}$.
We observe $\dagger $ is an antiautomorphism
of $\mathcal A$ such that $A^\dagger=A$
and $A^{*\dagger}=A^*$.
We have now shown there exists an
antiautomorphism
$\dagger $ of 
$\mathcal A$ such that
 $A^\dagger=A$
and $A^{*\dagger}=A^*$.
This antiautomorphism is unique since
$A,A^*$ together generate $\mathcal A$.
The map $X\rightarrow X^{\dagger \dagger}$ is
an isomorphism of $\K$-algebras from
$\mathcal A$ to itself.
This isomorphism is the identity
since
 $A^{\dagger \dagger}=A$,
 $A^{*\dagger \dagger}=A^*$,
 and since $A,A^*$ together generate 
$\mathcal A$.
\hfill $\Box $ \\

\begin{definition} 
\label{def:dag}
\rm
Let $A,A^*$ denote a Leonard pair in
 $\mathcal A$.
By the {\it antiautomorphism which corresponds to $A,A^*$}
 we mean the 
map $\dagger: {\mathcal A}\rightarrow {\mathcal A}$ from Theorem
\ref{thm:dagger}. Let 
$\Phi=(A;A^*;\lbrace E_i \rbrace_{i=0}^d; $ $  
\lbrace E^*_i\rbrace_{i=0}^d)$ denote a Leonard system in
 $\mathcal A$.
By the {\it antiautomorphism which corresponds to $\Phi$}
 we mean the
antiautomorphism which corresponds to the Leonard pair $A,A^*$. 
\end{definition}

\begin{lemma}
\label{cor:eistab}
Let $\Phi$ denote the Leonard system from Definition
\ref{eq:ourstartingpt} and let 
$\dagger $  denote the corresponding antiautomorphism.
Then the following (i), (ii) hold.
\begin{enumerate}
\item
Let $\mathcal D$ denote the subalgebra of $\cal A$ generated
by $A$. Then $X^\dagger=X$ for all $X \in \cal D$; in particular
$E_i^\dagger=E_i $ for $0 \leq i\leq d$. 
\item 
Let ${\mathcal D}^*$ denote the subalgebra of $\cal A$ generated
by $A^*$. Then $X^\dagger=X$ for all $X \in {\cal D}^*$; in particular
$E_i^{*\dagger}=E^*_i $ for $0 \leq i\leq d$. 
\end{enumerate}
\end{lemma}
\noindent {\it Proof:}
 (i) The sequence
$A^i$ $(0 \leq i \leq d)$ is a basis for
 the $\K$-vector space $\mathcal D$.
Observe $\dagger $ stabilizes $A^i$ for $0 \leq i \leq d$.
The result follows.
\\
\noindent (ii) Similar to the proof of (i) above. 
\hfill $\Box $ \\

\section{The scalars $a_i, x_i$}

\noindent In this section we introduce some scalars
that will help us describe Leonard systems.

\begin{definition}
\label{def:aixi}
\rm
Let $\Phi$ denote the Leonard system
from Definition 
\ref{eq:ourstartingpt}.
We define
\begin{eqnarray}
 a_i &=&\mbox{tr}(E^*_iA)\qquad    \quad   
  (0 \leq i \leq d),  \label{eq:aitr}
\\
x_i &=&\mbox{tr}(E^*_iAE^*_{i-1}A)     
\qquad  \quad   (1 \leq i \leq d),\label{eq:xitr} 
\end{eqnarray}
where $\mbox{tr}$ denotes trace.
For notational convenience we define $x_0=0$.
\end{definition}

\noindent 
We have a comment.

\begin{lemma}
\label{lem:bmat}
Let $\Phi$ denote the Leonard system
from Definition
\ref{eq:ourstartingpt} and let $V$ denote an irreducible
$\cal A$-module. For $0 \leq i \leq d$ let $v_i$ denote
a nonzero vector in $E^*_iV$ and observe $v_0, v_1, \ldots, v_d$
is a basis for $V$. Let $B$ denote the matrix in
$\hbox{Mat}_{d+1}(\K)$ that represents $A$ with respect to
$v_0, v_1, \ldots, v_d$. 
We observe $B$ is irreducible tridiagonal.
The following (i)--(iii) hold. 
\begin{enumerate}
\item
$B_{ii}=a_i$ $(0 \leq i \leq d)$.
\item
 $B_{i,i-1}B_{i-1,i}=x_i$  $(1 \leq i \leq d)$.
\item $x_i\not=0$ $(1 \leq i \leq d)$. 
\end{enumerate}
\end{lemma}

\noindent {\it Proof:}
(i), (ii) For $0 \leq i \leq d$ the matrix in
$\hbox{Mat}_{d+1}(\K)$ that represents $E^*_i$ with respect to
$v_0, v_1, \ldots, v_d$ has $ii$ entry
$1$ and all other entries 0. The result follows in
view of
Definition
\ref{def:aixi}.
\\
\noindent (iii) Immediate from (ii) and  since $B$ is irreducible.
\hfill $\Box $ \\

\begin{theorem}
\label{lem:monbasis}
Let $\Phi$ denote the Leonard system
from Definition
\ref{eq:ourstartingpt}. Let $V$ denote an irreducible
$\cal A$-module and  
let $v$ denote a nonzero vector in 
$E^*_0V$. Then
for $0 \leq i \leq d$ the vector
$E^*_iA^iv$ is nonzero and hence
a basis for $E^*_iV$. 
Moreover the sequence
\begin{eqnarray}
E^*_iA^iv \qquad (0 \leq i \leq d)
\label{eq:monbasis}
\end{eqnarray}
is a basis for $V$.
\end{theorem}
\noindent {\it Proof:}
We show $E^*_iA^iv\not=0$ for $0 \leq i \leq d$.
Let $i$ be given.
Setting $j=0$ in
Lemma
\ref{lem:absval}(iii)
we find
$E^*_iA^iE^*_0\not=0$.
Therefore
 $E^*_iA^iE^*_0V\not=0$.
The space $E^*_0V$  is spanned by
$v$
so $E^*_iA^iv\not=0$ as desired.
The remaining claims follow.
\hfill $\Box $ \\

\begin{theorem}
\label{lem:monbasis2}
Let $\Phi$ denote the Leonard system
from Definition
\ref{eq:ourstartingpt} and let
the scalars $a_i, x_i$ be as in
Definition
\ref{def:aixi}. Let
$V$ denote an irreducible
$\cal A$-module.
With respect to the basis for $V$ given
in 
(\ref{eq:monbasis})
 the matrix
that represents $A$ is equal to
\begin{equation}
\label{eq:monmat}
\left(
\begin{array}{ c c c c c c}
a_0 & x_1  &      &      &   &{\bf 0} \\
1 & a_1  &  x_2   &      &   &  \\
  & 1  &  \cdot    & \cdot  &   & \\
  &   & \cdot     & \cdot  & \cdot   & \\
  &   &           &  \cdot & \cdot & x_d \\
{\bf 0} &   &   &   & 1 & a_d  
\end{array}
\right) .
\end{equation} 
\end{theorem}

\noindent {\it Proof:}
With reference to
(\ref{eq:monbasis}) abbreviate
$v_i=E^*_iA^iv$ for $0 \leq i\leq d$.
Let $B$ denote the matrix in
$\hbox{Mat}_{d+1}(\K)$ that represents
$A$ with respect to
$v_0, v_1, \ldots, v_d$.
We show $B$ is equal to
(\ref{eq:monmat}). In view of
Lemma
\ref{lem:bmat} it suffices to show
$B_{i,i-1}=1$ for $1 \leq i \leq d$.
For $0 \leq i \leq d$ the 
matrix $B^i$ represents $A^i$ with
respect to 
$v_0, v_1, \ldots, v_d$; 
therefore 
$A^iv_0 = \sum_{j=0}^d(B^i)_{j0}v_j$.
Applying $E^*_i$ and using $v_0=v$
we find
$v_i = (B^i)_{i0}v_i$ so
$(B^i)_{i0}=1$.
By  Lemma
\ref{lem:step1}
we have $(B^i)_{i0}=
B_{i,i-1}
\cdots
B_{21}
B_{10}
$
so 
$
B_{i,i-1}
\cdots
B_{21}
B_{10}
=1$.
We now have
$
B_{i,i-1}
\cdots
B_{21}
B_{10}
=1$ for $1 \leq i \leq d$
so $B_{i,i-1}=1$ for $1 \leq i \leq d$. 
We now see $B$ is equal to
(\ref{eq:monmat}). 
\hfill $\Box $ \\

\begin{lemma}
\label{lem:xiprod}
Let $\Phi$ denote the Leonard system
from Definition
\ref{eq:ourstartingpt} and
let the scalars $a_i, x_i$ be as in 
Definition
\ref{def:aixi}.
Then the following (i)--(iii) hold.
\begin{enumerate}
\item
$E^*_iAE^*_i = a_iE^*_i \qquad (0 \leq i \leq d)$.
\item 
$E^*_iAE^*_{i-1}AE^*_i = x_iE^*_i \qquad \qquad (1 \leq i \leq d)$.
\item
$E^*_{i-1}AE^*_{i}AE^*_{i-1} = x_iE^*_{i-1} \qquad \qquad (1 \leq i \leq d)$.
\end{enumerate}
\end{lemma}

\noindent {\it Proof:}
(i) Setting $i=j$ and $r=0$ in Lemma
\ref{lem:absval}(iv) we find $E^*_i$ is a basis for
$E^*_i{\cal A}E^*_i$. By this and since
$E^*_i AE^*_i$ is contained in
$E^*_i{\cal A}E^*_i$ we find
there exists
 $\alpha_i \in \K$ such that
$E^*_i AE^*_i = \alpha_i E^*_i.$
Taking the trace of both sides and using
$\mbox{tr}(XY)=\mbox{tr}(YX)$, $\mbox{tr}(E^*_i)=1$ 
we
find $a_i = \alpha_i$.
\\
\noindent (ii) We mentioned above that
 $E^*_i$ is a basis for
$E^*_i{\cal A}E^*_i$. By this
and since 
$E^*_i AE^*_{i-1}AE^*_i$ is contained in
$E^*_i{\cal A}E^*_i$ we find
there exists
 $\beta_i \in \K$ such that
$E^*_i AE^*_{i-1}AE^*_i=\beta_i E^*_i$.
Taking the trace of both sides we 
find 
$x_i=\beta_i$.
\\
\noindent (iii)  Similar to the proof of (ii) above.
\hfill $\Box $ \\

\begin{lemma}
\label{lem:xrec}
Let $\Phi$ denote the Leonard system
from Definition  
\ref{eq:ourstartingpt} and 
let the scalars $ x_i$ be as in 
Definition
\ref{def:aixi}. Then the following (i), (ii) hold.
\begin{enumerate}
\item 
$ 
E^*_jA^{j-i}E^*_iA^{j-i}E^*_j = x_{i+1} x_{i+2} \cdots x_j E^*_j
\qquad (0 \leq i \leq j\leq d)$.
\item
$ 
E^*_iA^{j-i}E^*_jA^{j-i}E^*_i = x_{i+1} x_{i+2} \cdots x_j E^*_i
\qquad (0 \leq i \leq j\leq d)$.
\end{enumerate}
\end{lemma}
\noindent {\it Proof:}
(i) Evaluate the expression on the left using
Lemma
\ref{lem:absval}(ii), (iii) 
and Lemma
\ref{lem:xiprod}(ii).
\\
(ii)
Evaluate the expression on the left using
Lemma
\ref{lem:absval}(ii), (iii) 
and Lemma
\ref{lem:xiprod}(iii). 
\hfill $\Box $ \\

\section{The polynomials $p_i$}

\medskip
\noindent  
In this section we begin our discussion of 
polynomials. 
We will use
the following notation. 
Let $\lambda $ denote
an indeterminate. We let 
 $\K \lbrack \lambda \rbrack $ 
denote the $\K$-algebra consisting of all polynomials
in $\lambda $ that have coefficients in $\K$.  
For the rest of this paper all polynomials that we discuss
are assumed to lie in 
 $\K \lbrack \lambda \rbrack $.

\begin{definition}
\label{def:pi}
\rm
Let $\Phi$ denote the Leonard system
from Definition 
\ref{eq:ourstartingpt} and let
 the scalars $a_i, x_i$ be as in 
Definition
\ref{def:aixi}. We define
a sequence of 
polynomials $p_0, p_1, \ldots, p_{d+1}$ by
\begin{eqnarray}
p_0&=&1,
\\
\lambda p_i &=& p_{i+1}+a_ip_i+x_ip_{i-1}
\qquad (0 \leq i \leq d),
\label{eq:prec}
\end{eqnarray}
where $p_{-1}=0$.
We observe $p_i$ is monic with degree exactly $i$ for $0 \leq i \leq d+1$.
\end{definition}

\begin{lemma}
\label{lem:pareq}
Let $\Phi$ denote the Leonard system
from Definition 
\ref{eq:ourstartingpt} and let the polynomials $p_i$ be as in
Definition 
\ref{def:pi}.
Let $V$ denote
an irreducible $\mathcal A$-module and 
let $v$ denote a nonzero vector in $E^*_0V$.
Then
$p_i(A)v=E^*_iA^iv$ for $0 \leq i  \leq d$ and
$p_{d+1}(A)v=0$.
\end{lemma}
\noindent {\it Proof:}
We
abbreviate
$v_i=p_i(A)v$ for $0 \leq i\leq d+1$.
We define
$v'_i=E^*_iA^iv$ for $0 \leq i \leq d$
and $v'_{d+1}=0$.
We show $v_i=v'_i$ for $0 \leq i \leq d+1$.
From the construction
$v_0=v$ and $v'_0=v$ so $v_0=v'_0$.
From
(\ref{eq:prec}) we obtain
\begin{eqnarray}
Av_i = v_{i+1} + a_iv_i + x_iv_{i-1}
\qquad (0 \leq i \leq d)
\label{eq:rel}
\end{eqnarray}
where  $v_{-1}=0$.
From 
Theorem \ref{lem:monbasis2} we find
\begin{eqnarray}
Av'_i = v'_{i+1} + a_iv'_i + x_iv'_{i-1}
\qquad (0 \leq i \leq d)
\label{eq:rel2}
\end{eqnarray}
where  $v'_{-1}=0$.
Comparing
(\ref{eq:rel}), 
(\ref{eq:rel2}) and using $v_0=v'_0$
we find
$v_i=v'_i$ for $0 \leq i \leq d+1$. The result
follows.
\hfill $\Box $ \\

\noindent We mention a few consequences of Lemma
\ref{lem:pareq}.

\begin{theorem}
\label{thm:psend}
Let $\Phi$ denote the Leonard system
from Definition 
\ref{eq:ourstartingpt} and
let the 
polynomials $p_i$ be 
as in 
Definition
\ref{def:pi}. Let $V$ denote an irreducible 
$\mathcal A$-module.
 Then
$$
 p_i(A)E^*_0V = E^*_iV \quad (0 \leq i \leq d).
$$
\end{theorem}
\noindent {\it Proof:}
Let $v$ denote a nonzero vector in $E^*_0V$.
Then
$p_i(A)v=E^*_iA^iv$ by Lemma
\ref{lem:pareq}. 
Observe $v$ is a basis for $E^*_0V$.
By Theorem 
\ref{lem:monbasis} we find
$E^*_iA^iv$ is a basis for $E^*_iV$.
Combining these facts we find
 $p_i(A)E^*_0V = E^*_iV$.
\hfill $\Box $ \\

\begin{theorem}
\label{thm:pimon}
Let $\Phi$ denote the Leonard system
from Definition 
\ref{eq:ourstartingpt} and 
let the 
polynomials $p_i$ be 
as in 
Definition
\ref{def:pi}.
Then 
\begin{eqnarray} 
 p_i(A)E^*_0 = E^*_iA^iE^*_0 \qquad (0 \leq i \leq d).
\label{eq:pimon}
\end{eqnarray}
\end{theorem}
\noindent {\it Proof:}
Let the integer $i$ be given and abbreviate
$\Delta=
 p_i(A) -E^*_iA^i$.
We show
 $\Delta E^*_0=0$.
In order to do this we show
 $\Delta E^*_0V=0$,
 where $V$ denotes an irreducible $\mathcal A$-module.
Let $v$ denote a nonzero vector in $E^*_0V$ and recall
 $v$ is a basis for $E^*_0V$.
By Lemma
\ref{lem:pareq} we have
$\Delta v=0$ so
 $\Delta E^*_0V=0$.
Now 
 $\Delta E^*_0=0$ so
 $p_i(A)E^*_0 = E^*_iA^iE^*_0$.
\hfill $\Box $ \\

\begin{theorem}
\label{thm:pseq}
Let $\Phi$ denote the Leonard system
from Definition 
\ref{eq:ourstartingpt} and 
let the polynomial $p_{d+1}$ be as in
Definition
\ref{def:pi}. Then the following (i), (ii) hold.
\begin{enumerate}
\item 
 $p_{d+1}$ is both the minimal polynomial and the
 characteristic polynomial of
$A$.
\item
$
p_{d+1} = \prod_{i=0}^d (\lambda - \theta_i).
$
\end{enumerate}
\end{theorem}
\noindent {\it Proof:}
(i)
We first show $p_{d+1}$ is equal to the
minimal polynomial of $A$.
Recall
 $I, A, \ldots, A^d$ are linearly
independent and that
$p_{d+1}$ is monic with degree $d+1$.
We show
$p_{d+1}(A)=0$.
Let $V$ denote an irreducible $\mathcal A$-module.
Let $v$ denote a nonzero vector in $E^*_0V$
and recall $v$ is a basis for $E^*_0V$.
From Lemma
\ref{lem:pareq}
we find
$p_{d+1}(A)v=0$.
It follows
$p_{d+1}(A)E^*_0V=0$ so $p_{d+1}(A)E^*_0=0$.
Applying Lemma
\ref{lem:oneiszero}
(with $X=p_{d+1}(A)$ and $Y=I$)
we find $p_{d+1}(A)=0$.
We have now shown $p_{d+1}$ is the minimal polynomial
of $A$.
By definition 
the characteristic polynomial of $A$
is equal to 
$\mbox{det}(\lambda I-A)$.
This polynomial is monic with degree $d+1$ and
has $p_{d+1}$ as a factor; therefore
it is  equal to $p_{d+1}$.
\\
\noindent (ii) 
For $0 \leq i \leq d$ the scalar $\theta_i$
is an eigenvalue
of $A$ and therefore a root of the characteristic polynomial
of $A$.
\hfill $\Box $ \\

\begin{theorem}
\label{lem:eispoly}
Let $\Phi$ denote the Leonard system
from Definition 
\ref{eq:ourstartingpt} and
let the polynomials
 $p_i$ be 
as in 
Definition
\ref{def:pi}.
Let the 
scalars $x_i$ be as in
Definition
\ref{def:aixi}.
Then
\begin{eqnarray}
E^*_i = \frac{
p_i(A)E^*_0p_i(A)
}{
x_1x_2\cdots x_i
}
\qquad (0 \leq i \leq d).
\end{eqnarray}
\end{theorem}
\noindent {\it Proof:}
Let $\dagger:{\mathcal A}\rightarrow {\mathcal A}$ 
denote the antiautomorphism
which corresponds to 
$\Phi$.
From Theorem
\ref{thm:pimon}
we have 
$p_i(A)E^*_0
=
E^*_iA^iE^*_0
$.
Applying $\dagger$
we find
$E^*_0p_i(A)
=E^*_0A^iE^*_i
$.
From these comments  we find
\begin{eqnarray*}
p_i(A)E^*_0p_i(A) &=&
E^*_iA^iE^*_0
A^iE^*_i
\\
&=& x_1 x_2 \cdots x_iE^*_i
\end{eqnarray*}
in view of 
Lemma
\ref{lem:xrec}(i).
The result follows.
\hfill $\Box $ \\

\noindent 
We finish this section with a comment.

\begin{lemma}
\label{lem:aisum}
Let $\Phi$ denote the Leonard system
from Definition
\ref{eq:ourstartingpt} and
let the 
polynomials $p_i$ be 
as in 
Definition
\ref{def:pi}. Let the scalars $a_i$ be as in
Definition \ref{def:aixi}.
Then for $0 \leq i \leq d$ the coefficient of $\lambda^i$ in 
$p_{i+1}$ is equal to
$-\sum_{j=0}^i a_j $.
\end{lemma}

\noindent {\it Proof:}
Let $\alpha_i$ denote the coefficient of $\lambda^i$ in $p_{i+1}$.
Computing the coefficient of $\lambda^i$ in
(\ref{eq:prec}) we find $\alpha_{i-1}=\alpha_i+a_i$
for $0 \leq i \leq d$, where $\alpha_{-1}=0$. It follows
$\alpha_i = -\sum_{j=0}^i a_j$
for $0 \leq i \leq d$.
\hfill $\Box $ \\

\section{The scalars $\nu, m_i$}

\medskip
\noindent
In this section 
we introduce some more scalars that 
will help us describe Leonard systems.

\begin{definition} 
\label{def:mid}
\rm
Let $\Phi$ denote the Leonard system
from Definition  
\ref{eq:ourstartingpt}.
We define 
\begin{equation}
\label{eq:mid}
m_i =  \mbox{tr}(E_iE^*_0) \qquad \qquad (0 \leq i \leq d).
\end{equation}
\end{definition}

\begin{lemma}
\label{lem:mid}
Let $\Phi$ denote the Leonard system
from Definition
\ref{eq:ourstartingpt}.
Then (i)--(v) hold below.
\begin{enumerate}
\item $E_iE^*_0E_i = m_i E_i
\qquad (0 \leq i\leq d)$.
\item $E^*_0E_iE^*_0 = m_i E^*_0
\qquad (0 \leq i\leq d)$.
\item $m_i\not=0
\qquad (0 \leq i\leq d)$.
\item
$\sum_{i=0}^d m_i = 1$.
\item $m_0 = m^*_0$.
\end{enumerate}
\end{lemma}
\noindent {\it Proof:}
(i) Observe 
$E_i$ is a basis for $E_i{\mathcal A}E_i$.
By this and since
$E_iE^*_0E_i$ is contained in
 $E_i{\mathcal A}E_i$, there exists
 $\alpha_i \in \K$ 
such that 
 $E_iE^*_0E_i=\alpha_i E_i$. 
Taking the trace of both
 sides in this equation and using
$\hbox{tr}(XY)= \hbox{tr}(YX)$, $\mbox{tr}(E_i)=1$
  we find $\alpha_i=m_i$.
\\
\noindent (ii) Similar to the proof of (i).
\\
\noindent (iii) 
Observe $m_iE_i$ is equal to $E_iE^*_0E_i$ by part (i) above
and 
$E_iE^*_0E_i$ is nonzero by 
Corollary
\ref{cor:eibasis}. 
It follows $m_iE_i\not=0$ so $m_i\not=0$.
\\
\noindent (iv) 
Multiply each term
in the equation
$\sum_{i=0}^d E_i = I$
 on the right by $E^*_0$, and then take the trace.
Evaluate the result using
Definition
\ref{def:mid}.
\\
\noindent (v) 
The elements $E_0E^*_0$ and $E^*_0E_0$ have the same trace. 
\hfill $\Box $ \\

\begin{definition}
\label{def:n} 
\rm
 Let $\Phi$ denote the Leonard system
from Definition 
\ref{eq:ourstartingpt}. 
Recall $m_0=m^*_0$ by Lemma
\ref{lem:mid}(v); 
 we let $\nu$ denote the multiplicative inverse
of this common value. We observe
$\nu= \nu^*$. We 
emphasize
\begin{equation}
\mbox{tr}(E_0E^*_0) = \nu^{-1}.
\label{eq:sizedef}
\end{equation}
\end{definition}

\begin{lemma}
\label{lem:threeone}
Let $\Phi$ denote the Leonard 
system from Definition
\ref{eq:ourstartingpt}  and let 
the scalar $\nu$ be as in Definition
\ref{def:n}. Then the following (i), (ii) hold.
\begin{enumerate}
\item
$\nu E_0E^*_0E_0 = E_0$.
\item 
$\nu E^*_0E_0E^*_0 = E^*_0.$
\end{enumerate}
\end{lemma}
\noindent {\it Proof:}
(i) Set $i=0$ in Lemma
\ref{lem:mid}(i) and recall $m_0={\nu}^{-1}$.
\\
(ii) Set $i=0$ in Lemma
\ref{lem:mid}(ii) and recall $m_0={\nu}^{-1}$.
\hfill $\Box $ \\

\begin{theorem}
Let $\Phi$ denote the Leonard system
from Definition 
\ref{eq:ourstartingpt} and let the
polynomials $p_i$ be as in Definition \ref{def:pi}.
Let the scalars $\theta_i$ be as in
Definition \ref{def:evseq} and let the 
scalars $m_i$
be as in
Definition \ref{def:mid}.
Then 
\begin{eqnarray}
p_i(\theta_j) = m^{-1}_j \mbox{tr}(E_jE^*_iA^iE^*_0)
\qquad (0 \leq i,j\leq d).
\end{eqnarray}
\end{theorem}
\noindent {\it Proof:}
Using Theorem 
\ref{thm:pimon}
we find
\begin{eqnarray*}
\mbox{tr}(E_jE^*_iA^iE^*_0) &=& 
\mbox{tr}(E_jp_i(A)E^*_0)
\\
&=& p_i(\theta_j)\mbox{tr}(E_jE^*_0)
\\
&=& p_i(\theta_j)m_j.
\end{eqnarray*}
The result follows.
\hfill $\Box $ \\

\section{The standard basis}

In this section we 
 discuss the notion of a  {\it standard basis}.
We begin with a comment.

\begin{lemma}
\label{lem:stbasisp}
Let $\Phi$ denote the Leonard system
from Definition 
\ref{eq:ourstartingpt} and let $V$ denote an irreducible
$\mathcal A$-module.
Then 
\begin{eqnarray}
E^*_iV = E^*_iE_0V \qquad (0 \leq i \leq d).
\end{eqnarray}
\end{lemma}
\noindent {\it Proof:}
The space $E^*_iV$ has dimension 1 and contains
$E^*_iE_0V$. We show
$E^*_iE_0V\not=0$.
Applying
Corollary
\ref{cor:eibasis}
to $\Phi^*$ we find $E^*_iE_0\not=0$. It follows
 $E^*_iE_0V\not=0$.
We conclude
 $E^*_iV=E^*_iE_0V$.
\hfill $\Box $ \\

\begin{lemma}
\label{lem:stbasis}
Let $\Phi$ denote the Leonard system
from Definition 
\ref{eq:ourstartingpt} and
let
$V$ denote an irreducible $\mathcal A$-module.
Let 
$u$ denote a nonzero vector in  $E_0V$. Then
for $0 \leq i \leq d$ the vector $E^*_iu$ is nonzero
and hence a basis
for $E^*_iV$. Moreover the sequence
\begin{equation}
E^*_0u, E^*_1u, \ldots, E^*_du
\label{eq:stbasisint}
\end{equation}
is a basis for $V$.
\end{lemma}
\noindent {\it Proof:}
Let the integer $i$ be given.
We show $E^*_iu\not=0$.
Recall $E_0V$ has dimension 1 and
$u$ is a nonzero vector in 
 $E_0V$  
so $u$ spans  $E_0V$. Applying $E^*_i$ we find
$E^*_iu$ spans $E^*_iE_0V$. The space $E^*_iE_0V$ is nonzero by
Lemma
\ref{lem:stbasisp}
so 
$E^*_iu$ is nonzero.
The remaining assertions are clear.
\hfill $\Box $ \\

\begin{definition}
\label{def:stbasis}
\rm
Let $\Phi$ denote the Leonard system
from Definition 
\ref{eq:ourstartingpt} and let $V$ denote an irreducible
$\mathcal A$-module.
By a {\it $\Phi$-standard basis} for $V$,
we mean
a sequence 
$$
E^*_0u, E^*_1u, \ldots, E^*_du,
$$
where $u$ is a nonzero vector in $E_0V$.
\end{definition}

\noindent We give a few characterizations of the standard
basis.

\begin{lemma}
\label{lem:eggechar}
Let $\Phi$ denote the Leonard system from
Definition
\ref{eq:ourstartingpt}  and let $V$ denote an irreducible
$\mathcal A$-module.
Let $v_0, v_1, \ldots, v_d$ denote
a sequence of vectors in $V$, not all 0. Then this sequence is a 
$\Phi$-standard basis for $V$ if and only if both 
(i), (ii)  hold below.
\begin{enumerate}
\item $v_i \in E^*_iV$ for $0 \leq i \leq d$.
\item $\sum_{i=0}^d v_i \in E_0V$.
\end{enumerate}
\end{lemma} 
\noindent {\it Proof:}
To prove the lemma in one direction,
 assume 
$v_0, v_1, \ldots, v_d$ is a
$\Phi$-standard basis for $V$.
By Definition
\ref{def:stbasis}
there exists a nonzero $u \in E_0V$ such
that $v_i = E^*_iu$ for $0 \leq i \leq d$.
Apparently $v_i \in E^*_iV$ for $0 \leq i \leq d$ so (i) holds.
Let $I$ denote the identity element of $\mathcal A$ and
recall $I=\sum_{i=0}^d E^*_i$.
Applying this to $u$ we find
$u=\sum_{i=0}^d v_i $ and (ii) follows.
We have now proved the lemma in one direction. To prove the lemma
in the other direction, assume $v_0, v_1, \ldots, v_d$ satisfy
(i), (ii) above.
We define
$u=\sum_{i=0}^d v_i$ and observe $u\in E_0V$.
Using (i) 
we find $E^*_iv_j=\delta_{ij}v_j$ for
$0 \leq i,j\leq d$; it follows
$v_i = E^*_iu$ for $0 \leq i \leq d$.
Observe $u \not=0$  since at least one
of $v_0, v_1, \ldots, v_d$ is nonzero.
Now 
$v_0, v_1, \ldots, v_d$ is 
a 
$\Phi$-standard  basis for $V$ by Definition
\ref{def:stbasis}.
\hfill $\Box $ \\

\noindent  We recall some notation.
Let $d$ denote a nonnegative
integer and let $B$ denote a matrix in 
$\hbox{Mat}_{d+1}(\K)$. Let $\alpha $ denote a scalar in $\K$.
Then $B$ is said to have {\it constant row  sum $\alpha$} whenever
$B_{i0}+B_{i1}+\cdots + B_{id}=\alpha$ for $0 \leq i \leq d$.

\begin{lemma}
\label{lem:rowsum}
Let $\Phi$ denote the Leonard system from  Definition
\ref{eq:ourstartingpt} and let the scalars $\theta_i, \theta^*_i$
be as in Definition \ref{def:evseq}.
Let $V$ denote an irreducible $\mathcal A$-module and
let $v_0, v_1, \ldots, v_d$ denote a basis
for $V$.
Let $B$ (resp. $B^*$) denote the matrix in 
$\hbox{Mat}_{d+1}(\K)$
that represents
$A$ (resp. $A^*)$ with respect to 
this basis.
Then 
 $v_0, v_1, \ldots, v_d$ is a $\Phi$-standard basis
 for $V$ if and only if both (i), (ii) hold below.
\begin{enumerate}
\item $B$ has constant row sum $\theta_0$.
\item $B^*=\hbox{diag}(\theta^*_0, \theta^*_1, \ldots, \theta^*_d)$.
\end{enumerate}
\end{lemma}
\noindent {\it Proof:}
Observe $A \sum_{j=0}^d v_j = \sum_{i=0}^d v_i(B_{i0}+B_{i1}+\cdots B_{id})$.
Recall $E_0V$ is the eigenspace for $A$ and eigenvalue $\theta_0$.
Apparently 
$B$ has constant row sum $\theta_0$ if and only if
$\sum_{i=0}^d v_i \in E_0V$.
Recall that for $0 \leq i \leq d$,
$E^*_iV$ is the eigenspace for $A^*$ and eigenvalue $\theta^*_i$.
Apparently 
$B^*=\hbox{diag}(\theta^*_0, \theta^*_1, \ldots, \theta^*_d)$
if and only if $v_i \in E^*_iV$ for $0 \leq i \leq d$. 
The result follows in view of Lemma
\ref{lem:eggechar}.
\hfill $\Box $ \\

\begin{definition}
\label{def:flatcon}
\rm
Let $\Phi$ denote the Leonard system from 
Definition
\ref{eq:ourstartingpt}. We define
a map $\flat : {\mathcal A}\rightarrow 
\mbox{Mat}_{d+1}(\K)$ as follows.
Let 
$V$ denote an irreducible $\mathcal A$-module.
For all $X \in {\mathcal  A}$ we let $X^\flat $ denote
the matrix 
in 
$\mbox{Mat}_{d+1}(\K)$
that represents $X$ with respect to 
a $\Phi$-standard basis for $V$.
We observe $\flat : {\mathcal A} \rightarrow 
\mbox{Mat}_{d+1}(\K)$ is an isomorphism of
$\K$-algebras.
\end{definition}

\begin{lemma}
\label{lem:firstcom}
Let $\Phi$ denote the Leonard system from 
Definition
\ref{eq:ourstartingpt}  and let the scalars
$\theta_i, \theta^*_i$
be as in
Definition
\ref{def:evseq}.
Let the map $\flat : {\mathcal A}\rightarrow \mbox{Mat}_{d+1}(\K)$
be as in Definition
\ref{def:flatcon}.
Then (i)--(iii) hold below.
\begin{enumerate}
\item $A^\flat $ has constant row sum $\theta_0$.
\item $A^{*\flat}=\hbox{diag}(\theta^*_0,\theta^*_1, \ldots, \theta^*_d)$.
\item For $0 \leq i \leq d$ the matrix
$E_i^{*\flat}$ has $ii$ entry 1
and all other entries 0.
\end{enumerate}
\end{lemma}
\noindent {\it Proof:}
(i), (ii) Combine 
Lemma 
\ref{lem:rowsum} and
Definition
\ref{def:flatcon}.
\\
(iii) Immediate from
Lemma
\ref{lem:step0}(i).
\hfill $\Box $ \\

\noindent
Let $\Phi$ denote the Leonard system from 
Definition
\ref{eq:ourstartingpt}  and
let the map
$\flat : {\mathcal A}\rightarrow \mbox{Mat}_{d+1}(\K)$
be as in Definition
\ref{def:flatcon}.
Let $X$ denote an element of $\mathcal A$.
In 
Theorem
\ref{lem:firstcom3} below
we give
the entries  of $X^\flat$ in terms of the trace function.
To prepare for this we need a lemma.

\begin{lemma}
\label{lem:firstcom2}
Let $\Phi$ denote the Leonard system from 
Definition
\ref{eq:ourstartingpt}  and
let the map $\flat : {\mathcal A}\rightarrow \mbox{Mat}_{d+1}(\K)$
be as in Definition
\ref{def:flatcon}.
Then for $X \in \mathcal A$ the entries of $X^\flat$ satisfy
\begin{eqnarray}
E^*_iXE^*_jE_0 = (X^\flat)_{ij}E^*_iE_0
\qquad \qquad (0 \leq i,j\leq d).
\label{eq:xentrypre}
\end{eqnarray}
\end{lemma}
\noindent {\it Proof:}
Let the integers $i,j$ be given
and abbreviate
$\Delta=E^*_iXE^*_j-(X^\flat)_{ij}E^*_i$.
We show $\Delta E_0=0$. 
In order to do this we show $\Delta E_0V=0$, where
$V$ denotes an irreducible $\mathcal A$-module. Let $u$ denote
a nonzero vector in $E_0V$. 
By Definition
\ref{def:stbasis}
the sequence 
$E^*_0u, E^*_1u, \ldots, E^*_du$ 
is a 
 $\Phi$-standard basis 
for $V$.
Recall $X^\flat$ is the matrix in
$\mbox{Mat}_{d+1}(\K)$ that represents
$X$ with respect to
this basis.
Applying $X$ to $E^*_ju$ we find
$XE^*_ju = \sum_{r=0}^d 
(X^\flat)_{rj} E^*_ru$.
Applying 
$E^*_i$
we obtain
 $E^*_iXE^*_ju = 
(X^\flat)_{ij} E^*_iu$.
By this and since
$u$ spans $E_0V$ we find
$\Delta E_0V=0$.
Therefore $\Delta E_0=0$ and the result follows.
\hfill $\Box $ \\

\begin{theorem}
\label{lem:firstcom3}
Let $\Phi$ denote the Leonard system from 
Definition
\ref{eq:ourstartingpt}  and let the scalars
$m^*_i$ be as in
Definition
\ref{def:mid}.
Let the map $\flat : {\mathcal A}\rightarrow \mbox{Mat}_{d+1}(\K)$
be as in Definition
\ref{def:flatcon}. 
Then for $X\in {\mathcal A}$
 the entries of $X^\flat$ are given as follows.
\begin{eqnarray}
(X^\flat)_{ij}= m^{*-1}_i \mbox{tr}(E^*_iXE^*_jE_0)
\qquad \qquad (0 \leq i,j\leq d).
\label{eq:xentry}
\end{eqnarray}
\end{theorem}
\noindent {\it Proof:}
In 
equation
(\ref{eq:xentrypre}),
take the trace of both sides
and
observe $m^*_i=\mbox{tr}(E^*_iE_0)$ in
view of
Definition
\ref{def:mid}.
\hfill $\Box $ \\

\noindent Referring to Theorem
\ref{lem:firstcom3} we consider the case $X=E_0$.

\begin{lemma}
\label{lem:secondcom1}
Let $\Phi$ denote the Leonard system from 
Definition
\ref{eq:ourstartingpt}  and let the scalars
$m^*_i$ be as in
Definition
\ref{def:mid}.
Let the map $\flat : {\mathcal A}\rightarrow \mbox{Mat}_{d+1}(\K)$
be as in Definition
\ref{def:flatcon}.
Then for $0 \leq i,j\leq d$ the $ij$  entry
of $E_0^\flat$ is $m^*_j$.
\end{lemma}
\noindent {\it Proof:}
Set $X=E_0$ in (\ref{eq:xentry}).
Simplify the result using
$E_0E^*_jE_0=m^*_jE_0$ and
$m^*_i=\mbox{tr}(E^*_iE_0)$.
\hfill $\Box $ \\

\begin{theorem}
\label{thm:preimflat}
Let $\Phi$ denote the Leonard system from 
Definition
\ref{eq:ourstartingpt}, 
and let the map $\flat : {\mathcal A}\rightarrow \mbox{Mat}_{d+1}(\K)$
be as in Definition
\ref{def:flatcon}.
For $0 \leq i,j\leq d$ define
$\Psi_{ij}
=m^{*-1}_jE^*_iE_0E^*_j$,
where $m^*_j$ is from
Definition
\ref{def:mid}.
Then the matrix
$\Psi_{ij}^\flat$ 
has $ij$  entry
1 and all other entries 0.
\end{theorem}
\noindent {\it Proof:}
Immediate from
Lemma
\ref{lem:firstcom}(iii),
Lemma
\ref{lem:secondcom1},
and since $\flat $ is an isomorphism of $\K$-algebras.
\hfill $\Box $ \\

\section{The scalars $b_i, c_i$}

\noindent In this section we consider some scalars
that arise naturally in the context of the standard basis.

\begin{definition}
\label{def:sharpmp}
\rm
Let $\Phi$ denote the Leonard system
from Definition 
\ref{eq:ourstartingpt}
and let the map
$\flat :{\mathcal A}\rightarrow
\mbox{Mat}_{d+1}(\K)
$ be as in
Definition
\ref{def:flatcon}.
For $0 \leq i \leq d-1$ we let $b_i$ denote the
$i,i+1$ entry of $A^\flat$.
For
 $1 \leq i \leq d$ we let $c_i$ denote the
$i,i-1$ entry of $A^\flat$.
We observe 
\begin{equation}
\label{eq:matrepls}
A^\flat =  \left(
\begin{array}{ c c c c c c}
a_0 & b_0  &      &      &   &{\bf 0} \\
c_1 & a_1  &  b_1   &      &   &  \\
  & c_2  &  \cdot    & \cdot  &   & \\
  &   & \cdot     & \cdot  & \cdot   & \\
  &   &           &  \cdot & \cdot & b_{d-1} \\
{\bf 0} &   &   &   & c_d & a_d  
\end{array}
\right), 
\end{equation} 
where the $a_i$ are from 
Definition
\ref{def:aixi}.
For notational convenience we define $b_d=0$ and $c_0=0$.
\end{definition}

\begin{lemma}
\label{def:bici}
Let $\Phi$ denote the Leonard system
from Definition 
\ref{eq:ourstartingpt}
and let the scalars $b_i, c_i$ be as in
Definition
\ref{def:sharpmp}.
 Then with reference
to Definition
\ref{def:evseq} and Definition
\ref{def:aixi} the following (i), (ii) hold.
\begin{enumerate}
\item $b_{i-1}c_i=x_i \qquad (1 \leq i \leq d)$.
\item $c_i + a_i + b_i=\theta_0 \qquad (0 \leq i \leq d)$.
\end{enumerate}
\end{lemma}
\noindent {\it Proof:}
(i) Apply Lemma
\ref{lem:bmat}(ii) with $B=A^\flat$.
\\
(ii) Combine
(\ref{eq:matrepls})
and
Lemma
\ref{lem:firstcom}(i).
\hfill $\Box $ \\

\begin{lemma}
\label{def:bici2}
Let $\Phi$ denote the Leonard system
from Definition
\ref{eq:ourstartingpt} and let the scalars $b_i, c_i$
be as in
Definition
\ref{def:sharpmp}.
Let the polynomials $p_i$ be as in
Definition
\ref{def:pi}
and let the scalar $\theta_0$ be as in
Definition \ref{def:evseq}.
Then the following (i)--(iii) hold.
\begin{enumerate}
\item $b_i \not=0 \qquad (0 \leq i \leq d-1)$.
\item $c_i \not=0 \qquad (1 \leq i \leq d)$.
\item $b_0b_1 \cdots b_{i-1} = p_i(\theta_0) \qquad (0 \leq i \leq d+1)$.
\end{enumerate}
\end{lemma}
\noindent {\it Proof:}
(i), (ii) Immediate from 
Lemma
\ref{def:bici}(i) and
since each of $x_1, x_2, \ldots, x_d$ is nonzero.
\\
\noindent (iii) 
Assume $0 \leq i \leq d$; otherwise each side
is zero.
Let $\dagger:{\mathcal A}\rightarrow {\mathcal A}$
denote the antiautomorphism which
corresponds to $\Phi$.
Applying 
$\dagger $ to
both sides of
(\ref{eq:pimon})
we get
$E^*_0p_i(A)=
E^*_0A^{i}E^*_i$.
We may now argue
\begin{eqnarray*}
b_0b_1\ldots b_{i-1}
&=&
(A^{i\flat})_{0i}
\qquad \qquad  
\qquad \qquad  
(\hbox{by (\ref{eq:matrepls})})
\\
&=& m^{*-1}_0\mbox{tr}(E^*_0A^iE^*_iE_0)
\qquad \qquad  
(\hbox{by 
Theorem
\ref{lem:firstcom3}})
\\
&=& m^{*-1}_0
\mbox{tr}(E^*_0p_i(A)E_0)
\\
&=& m^{*-1}_0p_i(\theta_0)
\mbox{tr}(E^*_0E_0)
\\
&=& p_i(\theta_0)
\qquad \qquad 
\qquad \qquad 
(\mbox{by 
Definition \ref{def:mid}}).
\end{eqnarray*}
\hfill $\Box $ \\

\begin{theorem}
\label{def:bici3}
Let $\Phi$ denote the Leonard system
from Definition 
\ref{eq:ourstartingpt} and let the polynomials
$p_i$ be as in
Definition
\ref{def:pi}. Let the scalar $\theta_0$ be as in
Definition \ref{def:evseq}.
Then 
$p_i(\theta_0)\not=0$ for 
$0 \leq i \leq d$.
Let the scalars $b_i, c_i$ be
as in
Definition
\ref{def:sharpmp}.
Then 
\begin{eqnarray}
\label{eqbi}
b_i&=& \frac{p_{i+1}(\theta_0)}
{p_i(\theta_0)} 
\qquad  \qquad 
(0 \leq i \leq d)
\end{eqnarray}
and 
\begin{eqnarray}
\label{eqci}
c_i&=& \frac{x_ip_{i-1}(\theta_0)}
{p_i(\theta_0)} 
\qquad  \qquad 
(1 \leq i \leq d).
\end{eqnarray}
\end{theorem}
\noindent {\it Proof:}
Observe $p_i( \theta_0)\not=0$ for $0 \leq i \leq d$
by Lemma
\ref{def:bici2}(i), (iii).
Line 
(\ref{eqbi}) is immediate from
Lemma
\ref{def:bici2}(iii).
To get 
(\ref{eqci}) combine
(\ref{eqbi}) and
Lemma \ref{def:bici}(i).
\hfill $\Box $ \\

\begin{lemma}
\label{def:bici4a}
Let $\Phi$ denote the Leonard system
from Definition 
\ref{eq:ourstartingpt} and let the scalars $b_i, c_i$
be as in Definition
\ref{def:sharpmp}.
Then the following (i), (ii) hold.
\begin{enumerate}
\item 
$
E^*_{i}AE^*_{i+1}E_0 = b_iE^*_iE_0 \qquad 
(0 \leq i \leq d-1).
$
\item
$
E^*_{i}AE^*_{i-1}E_0
=c_iE^*_iE_0
\qquad 
(1 \leq i \leq d).
$
\end{enumerate}
\end{lemma}
\noindent {\it Proof:}
(i) 
This is (\ref{eq:xentrypre})
with
$X=A$ and $j=i+1$.
\\
\noindent (ii) 
This is
(\ref{eq:xentrypre})
with
$X=A$ and $j=i-1$.
\hfill $\Box $ \\

\begin{theorem}
\label{def:bici4}
Let $\Phi$ denote the Leonard system
from Definition 
\ref{eq:ourstartingpt} and let the scalars $b_i, c_i$
be as in Definition
\ref{def:sharpmp}. Let the scalars $m^*_i$ be as in
Definition
\ref{def:mid}.
Then the following (i), (ii) hold.
\begin{enumerate}
\item 
$
b_i = m^{*-1}_i tr(E^*_{i}AE^*_{i+1}E_0) \qquad 
(0 \leq i \leq d-1).
$
\item
$
c_i = m^{*-1}_i tr(E^*_{i}AE^*_{i-1}E_0) \qquad 
(1 \leq i \leq d).
$
\end{enumerate}
\end{theorem}
\noindent {\it Proof:}
(i) 
This is
(\ref{eq:xentry})
with $X=A$ and $j=i+1$.
\\
\noindent (ii)
This is
(\ref{eq:xentry})
with $X=A$ and $j=i-1$.
\hfill $\Box $ \\

\noindent We finish this section with a comment.

\begin{theorem}
\label{lem:longeq}
Let $\Phi$ denote the Leonard system
from Definition 
\ref{eq:ourstartingpt} and let the scalars $ c_i$ be as in
Definition
\ref{def:sharpmp}.
Let the scalars $\theta_i$
be as in Definition \ref{def:evseq}
and let the scalar $\nu$ be as in
Definition \ref{def:n}.
 Then
\begin{equation}
\label{frame}
(\theta_0-\theta_1)
(\theta_0-\theta_2) \cdots
(\theta_0-\theta_d)
= 
\nu c_1c_2\cdots c_d.
\end{equation}
\end{theorem}
\noindent {\it Proof:}
Let $\delta$ denote the expression on the left-hand side
of 
(\ref{frame}). Setting $i=0$ in
(\ref{eq:primiddef}) we find
$\delta E_0=\prod_{j=1}^d(A-\theta_jI)$.
We multiply both sides of this equation on the left by
$E^*_d$ and on the right by $E^*_0$. We evaluate
the resulting equation using
Lemma
\ref{lem:absval}(i)
to obtain
$\delta E^*_dE_0E^*_0=E^*_dA^dE^*_0$.
We multiply both sides of this equation on the right by $E_0$ to obtain
\begin{eqnarray}
\label{eq:fourterm}
\delta E^*_dE_0E^*_0E_0=E^*_dA^dE^*_0E_0.
\end{eqnarray}
We evaluate each side of 
(\ref{eq:fourterm}).
The left-hand side of 
(\ref{eq:fourterm}) is equal to
$\delta \nu^{-1} E^*_dE_0$ in view of
Lemma
\ref{lem:threeone}(i). 
We now consider the right-hand side of 
(\ref{eq:fourterm}).  Observe
$E^*_dA^dE^*_0=E^*_dAE^*_{d-1}A$ 
$\cdots  
 E^*_1AE^*_0$
by (\ref{eq:rightside}).
Evaluating the right-hand side of
(\ref{eq:fourterm}) using this
 and Theorem
\ref{def:bici4a}(ii)
we find 
it is equal to 
$c_1c_2\cdots c_d E^*_dE_0$.
From our above comments
we find
$\delta \nu^{-1} E^*_dE_0=
c_1c_2\cdots c_d E^*_dE_0$.
Observe $E^*_dE_0\not=0$ by
Lemma
\ref{lem:stbasisp}
so 
$\delta \nu^{-1}=
c_1c_2\cdots c_d$.
The result follows.
\hfill $\Box $ \\

\section{The scalars $k_i$}

\noindent In this section we consider some scalars that are closely
related to the scalars from Definition
\ref{def:mid}.

\begin{definition}
\label{def:ki1}
\rm
Let $\Phi$ denote the Leonard system
from Definition
\ref{eq:ourstartingpt}.
We define
\begin{eqnarray}
\label{eq:ki1}
k_i = m^*_i \nu \qquad   (0 \leq i \leq d),
\end{eqnarray}
where the $m^*_i$ are from
Definition
\ref{def:mid} and $\nu$ is from
Definition
\ref{def:n}. 
\end{definition}

\begin{lemma}
\label{def:ki2}
Let $\Phi$ denote the Leonard system
from
Definition \ref{eq:ourstartingpt} and let the scalars $k_i$
be as in
Definition
\ref{def:ki1}. Then
(i) $k_0=1$;
(ii) $k_i\not=0$ for  $0 \leq i \leq d$;
(iii) $\sum_{i=0}^d k_i  = \nu$.

\end{lemma}
\noindent {\it Proof:}
(i) Set $i=0$ in 
(\ref{eq:ki1}) and recall $m^*_0=\nu^{-1}$.
\\
(ii)  Applying Lemma
\ref{lem:mid}(iii) to $\Phi^*$ we find
$m^*_i\not=0$ for $0 \leq i \leq d$.
We have $\nu\not=0$ by Definition
\ref{def:n}. 
The result follows in view of
(\ref{eq:ki1}).
\\
\noindent (iii)
  Applying Lemma
\ref{lem:mid}(iv) to $\Phi^*$ we find
$\sum_{i=0}^d m^*_i=1$. The result follows
in view of 
(\ref{eq:ki1}).
\hfill $\Box $ \\

\begin{lemma}
\label{def:ki3}
Let $\Phi$ denote the Leonard system
from
Definition \ref{eq:ourstartingpt} and let the scalars $k_i$
be as in
Definition
\ref{def:ki1}.
Then with reference to Definition
\ref{def:evseq}, Definition 
\ref{def:aixi}, and Definition
\ref{def:pi},
\begin{eqnarray}
\label{eq:kipi}
k_i = \frac
{
p_i(\theta_0)^2
}
{
x_1 x_2 \cdots x_i
}
\qquad (0 \leq i \leq d).
\end{eqnarray}
\end{lemma}
\noindent {\it Proof:}
We show that each side of 
(\ref{eq:kipi}) is equal to
$\nu {\rm tr}(E^*_iE_0)$.
Using 
(\ref{eq:mid}) and
(\ref{eq:ki1}) we find
$\nu {\rm tr}(E^*_iE_0)$ is equal
to the left-hand side of
(\ref{eq:kipi}).
Using Theorem \ref{lem:eispoly} we 
find 
$\nu {\rm tr}(E^*_iE_0)$ is equal
to the right-hand side of
(\ref{eq:kipi}).
\hfill $\Box $ \\

\begin{theorem}
\label{def:ki4}
Let $\Phi$ denote the Leonard system
from
Definition \ref{eq:ourstartingpt} and let the scalars $k_i$
be as in
Definition
\ref{def:ki1}. Let the scalars $b_i, c_i$ be as in
Definition 
\ref{def:sharpmp}.
Then
\begin{eqnarray}
k_i = \frac
{
b_0b_1\cdots b_{i-1}
}
{
c_1 c_2 \cdots c_i
}
\qquad (0 \leq i \leq d).
\label{eq:kibc}
\end{eqnarray}
\end{theorem}
\noindent {\it Proof:}
Evaluate the expression on the right in
(\ref{eq:kipi}) using
Lemma \ref{def:bici}(i) and 
Lemma \ref{def:bici2}(iii).
\hfill $\Box $ \\

\section{The polynomials  $v_i$}

\noindent Let $\Phi$ denote the Leonard system from Definition
\ref{eq:ourstartingpt} and let the polynomials $p_i$ be
as in
Definition \ref{def:pi}.
The  $p_i$ have two normalizations
 of interest; we call these the
$u_i$ and the $v_i$. In this section we discuss the
$v_i$.
In the next section we will discuss the $u_i$.

\begin{definition}
\label{def:vi1}
\rm
Let $\Phi$ denote the Leonard system
from Definition 
\ref{eq:ourstartingpt} and let the polynomials $p_i$ be as in
 Definition
\ref{def:pi}.
For $0 \leq i \leq d$ we define the polynomial $v_i $ 
by 
\begin{equation}
\label{vicl}
v_i = \frac{p_i}{c_1c_2\cdots c_i},
\end{equation}
where the 
 $c_j$ are
from Definition
\ref{def:sharpmp}.
We observe  $v_0=1$.
\end{definition}

\begin{lemma}
\label{def:vi1a}
Let $\Phi$ denote the Leonard system
from Definition 
\ref{eq:ourstartingpt} and let the polynomials $v_i$ be as in
 Definition 
\ref{def:vi1}. Let the 
scalar $\theta_0$ be as in Definition \ref{def:evseq} and let the
scalars $k_i$ be as in
Definition
\ref{def:ki1}.
Then
\begin{eqnarray}
\label{eq:vnorm}
 v_i(\theta_0)= k_i
\qquad ( 0 \leq i \leq d).
\end{eqnarray}
\end{lemma}
\noindent {\it Proof:}
Use Lemma
\ref{def:bici2}(iii),
Theorem
\ref{def:ki4}, and
(\ref{vicl}).
\hfill $\Box $ \\

\begin{lemma}
\label{def:vi2}
Let $\Phi$ denote the Leonard system
from Definition 
\ref{eq:ourstartingpt} 
and let the polynomials $v_i$ be as in
 Definition 
\ref{def:vi1}. Let the scalars $a_i, b_i, c_i$ be as in
Definition
\ref{def:aixi} and Definition
\ref{def:sharpmp}.
Then
\begin{eqnarray}
\lambda v_i = c_{i+1}v_{i+1}
+ a_i v_i 
+ b_{i-1} v_{i-1}
\qquad (0 \leq i \leq d-1),
\end{eqnarray}
where $b_{-1}=0$ and $v_{-1}=0$. Moroever
\begin{eqnarray}
\lambda v_d 
- a_d v_d 
-b_{d-1}v_{d-1} = (c_1 c_2 \cdots c_d)^{-1}p_{d+1}.
\end{eqnarray}
\end{lemma}
\noindent {\it Proof:}
In 
(\ref{eq:prec}), 
divide both sides  by $c_1c_2 \cdots c_i$.
Evaluate the result using
Lemma \ref{def:bici}(i)
and 
(\ref{vicl}).
\hfill $\Box $ \\

\begin{theorem}
\label{def:vi3}
Let $\Phi$ denote the Leonard system
from Definition
\ref{eq:ourstartingpt}
and let the polynomials $v_i$ be as in 
 Definition 
\ref{def:vi1}.
Let $V$ denote an irreducible $\mathcal A$-module and let
$u$ denote a nonzero vector in $E_0V$.
Then
\begin{equation}
\label{eq:viau}
v_i(A)E^*_0u = E^*_iu \qquad \qquad (0 \leq i \leq d).
\end{equation}
\end{theorem}
\noindent {\it Proof:}
For $0 \leq i \leq d$ we define
$w_i=v_i(A)E^*_0u$ 
and $w'_i=E^*_iu$.
We show $w_i=w'_i$.
Each of $w_0$,  $w'_0$ is equal to
$E^*_0u$ so $w_0=w'_0$.
Using
Lemma \ref{def:vi2}
we obtain
\begin{eqnarray}
Aw_i = c_{i+1}w_{i+1} + a_iw_i + b_{i-1}w_{i-1}
\qquad (0 \leq i \leq d-1)
\label{eq:releta}
\end{eqnarray}
where  $w_{-1}=0$ and $b_{-1}=0$.
By
Definition 
\ref{def:flatcon},
Definition
\ref{def:sharpmp}, and since
$w'_0, w'_1, \ldots, w'_d$
is a $\Phi$-standard basis,
\begin{eqnarray}
Aw'_i = c_{i+1}w'_{i+1} + a_iw'_i + b_{i-1}w'_{i-1}
\qquad (0 \leq i \leq d-1)
\label{eq:rel2eta}
\end{eqnarray}
where  $w'_{-1}=0$.
Comparing
(\ref{eq:releta}), 
(\ref{eq:rel2eta}) and using $w_0=w'_0$
we find
$w_i=w'_i$ for $0 \leq i \leq d$.
The result follows.
\hfill $\Box $ \\

\noindent We finish this section with a comment.

\begin{lemma}
\label{def:vi4}
Let $\Phi$ denote the Leonard system
from Definition
\ref{eq:ourstartingpt}
and let the polynomials $v_i$ be as in 
 Definition 
\ref{def:vi1}. Let the scalar $\nu$ be as in 
Definition \ref{def:n}.
Then the following (i), (ii) hold.
\begin{enumerate}
\item
$v_i(A)E^*_0E_0 = E^*_iE_0 \qquad (0 \leq i \leq d).$
\item
$v_i(A)E^*_0 = \nu E^*_iE_0E^*_0 \qquad (0 \leq i \leq d).$
\end{enumerate}
\end{lemma}
\noindent {\it Proof:}
(i) Let the integer $i$ be given and abbreviate
$
\Delta=
v_i(A)E^*_0-E^*_i
$.
We show $\Delta E_0=0$.
In order to to do this we show  
$\Delta E_0V=0$, where $V$ denotes an irreducible 
$\mathcal A$-module.
Let $u$ denote a nonzero vector in $E_0V$ and 
recall
$u$ spans $E_0V$.
Observe
$\Delta u=0$ by 
Theorem \ref{def:vi3} so
$\Delta E_0V=0$. Now $\Delta E_0=0$ so
$v_i(A)E^*_0E_0 = E^*_iE_0$.
\\
\noindent (ii) In the equation of (i) above, multiply  both sides
on the right by $E^*_0$ and simplify the result using
Lemma \ref{lem:threeone}(ii).
\hfill $\Box $ \\

\section{The polynomials $u_i$}

\noindent Let $\Phi$ denote the Leonard system from Definition
\ref{eq:ourstartingpt} and let the polynomials $p_i$ be
as in
Definition \ref{def:pi}. In the previous section we gave
a 
 normalization of the $p_i$ that we called the $v_i$.
 In this section we give a second normalization for
 the $p_i$ that we call the $u_i$.

\begin{definition}
\label{def:ui1}
\rm
Let $\Phi$ denote the Leonard system
from 
Definition 
\ref{eq:ourstartingpt} and let the polynomials $p_i$ be
as in Definition
\ref{def:pi}.
For $0 \leq i \leq d$ we define 
the polynomial $u_i$ by
\begin{equation}
\label{uicl}
u_i = \frac{p_i}{p_i(\theta_0)},
\end{equation}
where $\theta_0$ is from
Definition
\ref{def:evseq}.
We observe $u_0=1$. Moreover 
\begin{eqnarray}
\label{eq:unorm}
u_i(\theta_0)= 1 
\qquad ( 0 \leq i \leq d).
\end{eqnarray}
\end{definition}

\begin{lemma}
\label{def:ui2}
Let $\Phi$ denote the Leonard system
from 
Definition 
\ref{eq:ourstartingpt}
and let the polynomials
$u_i$ be as in
Definition 
\ref{def:ui1}.
Let the scalars $a_i, b_i, c_i$ be as in
Definition
\ref{def:aixi} and Definition
\ref{def:sharpmp}.
Then  
\begin{eqnarray}
\lambda u_i = b_iu_{i+1}
+ a_i u_i 
+ c_i u_{i-1}
\qquad (0 \leq i \leq d-1),
\label{eq:uirec1}
\end{eqnarray}
where $u_{-1}=0$.
Moreover
\begin{eqnarray}
\lambda u_d 
- c_d u_{d-1} 
-a_du_d = p_d(\theta_0)^{-1}p_{d+1},
\label{eq:uirec2}
\end{eqnarray}
where $\theta_0$ is from 
Definition \ref{def:evseq}.
\end{lemma}
\noindent {\it Proof:}
In 
(\ref{eq:prec}), 
divide both sides  by $p_i(\theta_0)$
and evaluate the result using
Lemma \ref{def:bici}(i),
(\ref{eqbi}),
and (\ref{uicl}).
\hfill $\Box $ \\

\noindent The above 3-term recurrence is often expressed
as follows.

\begin{corollary}
\label{cor:3term}
Let $\Phi$ denote the Leonard system
from 
Definition 
\ref{eq:ourstartingpt}
and let the polynomials
$u_i$ be as in
Definition 
\ref{def:ui1}. Let the scalars $\theta_i$
be as in Definition \ref{def:evseq}.
Then  for $0 \leq i,j\leq d$ we have 
\begin{eqnarray}
\theta_j u_i(\theta_j) = b_iu_{i+1}(\theta_j)
+ a_i u_i(\theta_j)
+ c_i u_{i-1}(\theta_j),
\label{eq:3termu}
\end{eqnarray}
where $u_{-1}=0$ and $u_{d+1}=0$.
\end{corollary}
\noindent {\it Proof:}
Apply Lemma
\ref{def:ui2}
(with $\lambda = \theta_j$)
and observe $p_{d+1}(\theta_j)=0$ by Theorem
\ref{thm:pseq}(ii).
\hfill $\Box $ \\

\begin{lemma}
\label{def:ui4}
Let $\Phi$ denote the Leonard system
from 
Definition 
\ref{eq:ourstartingpt}.
Let the polynomials $u_i, v_i$ be as in
Definition
\ref{def:ui1} and 
Definition
\ref{def:vi1} respectively.
Then 
\begin{equation}
\label{eq:uvk}
v_i = k_i u_i  \qquad (0 \leq i \leq d),
\end{equation}
where the $k_i$ are from
Definition \ref{def:ki1}.
\end{lemma}
\noindent {\it Proof:}
Compare 
(\ref{vicl}) and
(\ref{uicl}) in light of
Lemma \ref{def:bici2}(iii) and
Theorem
\ref{def:ki4}.
\hfill $\Box $ \\

\noindent Let $\Phi$ denote the 
Leonard system
from Definition
\ref{eq:ourstartingpt} and 
let the polynomials
$u_i$ be as in
Definition
\ref{def:ui1}.
Let $\theta_0, \theta_1, \ldots, \theta_d$ denote
the eigenvalue sequence of $\Phi$.
Our next goal is
to compute the $u_i(\theta_j)$ in terms
of the trace function.
To prepare for this we give a lemma.

\begin{lemma}
\label{def:vi5pre}
Let $\Phi$ denote the Leonard system
from 
Definition 
\ref{eq:ourstartingpt} and 
let the polynomials $v_i$ be as in
Definition
\ref{def:vi1}. Let the scalars $\theta_i$
 be as in
Definition
\ref{def:evseq} and let the scalars $m_i$ be as in
Definition
\ref{def:mid}.
Then 
\begin{eqnarray}
E_0E^*_iE_jE^*_0 = v_i(\theta_j)m_j E_0E^*_0
 \qquad (0 \leq i,j \leq d).
\label{eq:uijpre}
\end{eqnarray}
\end{lemma}
\noindent {\it Proof:}
Let $\dagger :{\mathcal A}\rightarrow {\mathcal A}$
denote the antiautomorphism 
which corresponds to $\Phi$.
Applying $\dagger $ to
the equation in Lemma
\ref{def:vi4}(i) we find
$
E_0E^*_0v_i(A)=
E_0E^*_i$.
Using  this and 
Lemma
\ref{lem:mid}(ii) 
we find
\begin{eqnarray*}
E_0E^*_iE_jE^*_0 &=& 
E_0E^*_0v_i(A)E_jE^*_0\\
&=& v_i(\theta_j)E_0E^*_0E_jE^*_0\\
&=& v_i(\theta_j)m_jE_0E^*_0.
\end{eqnarray*}
\hfill $\Box $ \\

\begin{theorem}
\label{def:vi5}
Let $\Phi$ denote the Leonard system
from 
Definition 
\ref{eq:ourstartingpt} and 
let the polynomials $u_i$ be as in
Definition
\ref{def:ui1}. 
Let the scalars $\theta_i$ be as in
Definition
\ref{def:evseq} and let the scalars
 $m_i, m^*_i$ be as in
Definition
\ref{def:mid}.
Then 
\begin{eqnarray}
u_i(\theta_j) = m^{*-1}_i m^{-1}_j \mbox{tr}(E_0E^*_iE_jE^*_0)
 \qquad (0 \leq i,j \leq d).
\label{eq:uiform}
\end{eqnarray}
\end{theorem}
\noindent {\it Proof:}
In (\ref{eq:uijpre}), 
take the trace of both sides
 and simplify the result using
(\ref{eq:sizedef}),
(\ref{eq:ki1}),
(\ref{eq:uvk}).
\hfill $\Box $ \\

\begin{theorem}
\label{lem:uidual}
Let $\Phi$ denote the Leonard system
from 
Definition 
\ref{eq:ourstartingpt}. Let the polynomials
$u_i$ be as in Definition
\ref{def:ui1} and recall the $u^*_i$ are the 
corresponding polynomials for $\Phi^*$.
Let the scalars $\theta_i, \theta^*_i$ be as in
Definition \ref{def:evseq}.
Then
\begin{eqnarray}
u_i(\theta_j) = u^*_j(\theta^*_i) 
\qquad (0 \leq i,j\leq d).
\label{eq:uijus}
\end{eqnarray}
\end{theorem}
\noindent {\it Proof:}
Applying Theorem 
\ref{def:vi5} to $\Phi^*$ we find
\begin{eqnarray}
u^*_i(\theta^*_j) = m^{-1}_i m^{*-1}_j \mbox{tr}(E^*_0E_iE^*_jE_0).
\label{eq:uisbeg}
\end{eqnarray}
Interchanging the roles of $i,j$ in 
(\ref{eq:uisbeg})
 we obtain
\begin{eqnarray}
u^*_j(\theta^*_i) = m^{*-1}_i m^{-1}_j \mbox{tr}(E^*_0E_jE^*_iE_0).
\label{eq:uiforms}
\end{eqnarray}
Let $\dagger :{\mathcal A}\rightarrow {\mathcal A} $ 
denote the antiautomorphism 
which corresponds to $\Phi$.
Observe
\begin{eqnarray}
\label{eq:daggeref}
(E_0E^*_iE_jE^*_0)^\dagger = E^*_0E_jE^*_iE_0
\end{eqnarray}
in view of  Lemma
\ref{cor:eistab}.
The trace function is invariant under $\dagger$ so
\begin{eqnarray}
\label{eq:trdagger}
\mbox{tr}(E_0E^*_iE_jE^*_0)= \mbox{tr}(E^*_0E_jE^*_iE_0).
\end{eqnarray}
 Combining 
(\ref{eq:uiform}), 
(\ref{eq:uiforms}),
(\ref{eq:trdagger})
we obtain
(\ref{eq:uijus}).
\hfill $\Box $ \\

\noindent In the following two theorems we show how 
(\ref{eq:uijus}) looks in terms of the polynomials
$v_i$ and $p_i$.

\begin{theorem}
\label{lem:vidual}
Let $\Phi$ denote the Leonard system
from 
Definition 
\ref{eq:ourstartingpt}.
With reference to 
 Definition
\ref{def:evseq},
Definition \ref{def:notconv}
and Definition
\ref{def:vi1},
\begin{eqnarray}
v_i(\theta_j)/k_i = v^*_j(\theta^*_i)/k^*_j 
\qquad (0 \leq i,j\leq d).
\label{eq:vijus}
\end{eqnarray}
\end{theorem}
\noindent {\it Proof:}
Evaluate (\ref{eq:uijus})
using
Lemma
\ref{def:ui4}.
\hfill $\Box $ \\

\begin{theorem}
\label{lem:pidual}
Let $\Phi$ denote the Leonard system
from 
Definition 
\ref{eq:ourstartingpt}.
With reference to 
 Definition
\ref{def:evseq},
Definition \ref{def:notconv},
and Definition
\ref{def:pi},
\begin{eqnarray}
\frac{p_i(\theta_j)}{p_i(\theta_0)}
= 
\frac{p^*_j(\theta^*_i)}{p^*_j(\theta^*_0)}
\qquad (0 \leq i,j\leq d).
\label{eq:pijus}
\end{eqnarray}
\end{theorem}
\noindent {\it Proof:}
Evaluate (\ref{eq:uijus}) 
using
Definition \ref{def:ui1}.
\hfill $\Box $ \\

\noindent The equations
(\ref{eq:uijus}),
(\ref{eq:vijus}),
(\ref{eq:pijus})
are often referred to as {\it Askey-Wilson duality}.

\medskip
\noindent 
We finish this section with a few comments.

\begin{lemma}
\label{lem:diffeq}
Let $\Phi$ denote the Leonard system
from 
Definition 
\ref{eq:ourstartingpt} and
let the polynomials $u_i$ be as in
 Definition 
\ref{def:ui1}.
Then for $0 \leq i,j\leq d$ we have
\begin{eqnarray}
\label{eq:diffeq}
\theta^*_i u_i(\theta_j) = 
b^*_ju_i(\theta_{j+1})
+
a^*_ju_i(\theta_j)
+
c^*_ju_i(\theta_{j-1}),
\end{eqnarray}
where $\theta_{-1}$, $\theta_{d+1}$ denote indeterminates.
\end{lemma}
\noindent {\it Proof:}
Apply 
Corollary \ref{cor:3term}
to $\Phi^*$ and evaluate the result using
Theorem
\ref{lem:uidual}.
\hfill $\Box $ \\

\noindent We refer to 
(\ref{eq:diffeq}) as the {\it difference equation}
satisfied by the $u_i$.

\begin{lemma}
\label{lem:thsform}
Let $\Phi$ denote the Leonard system
from 
Definition 
\ref{eq:ourstartingpt} and assume $d\geq 1$.
Let the polynomials $u_i$ be as in
 Definition 
\ref{def:ui1}.
Then
\begin{eqnarray}
u_i(\theta_1) = \frac{\theta^*_i-a^*_0}{\theta^*_0-a^*_0}
\qquad (0 \leq i\leq d).
\label{eq:thsform}
\end{eqnarray}
\end{lemma}
\noindent {\it Proof:}
Setting $j=1$ in 
(\ref{eq:uijus}) we find
$u_i(\theta_1)=u^*_1(\theta^*_i).
$
Applying 
(\ref{uicl}) to $\Phi^*$ we find
$u^*_1=p^*_1/p^*_1(\theta^*_0)$.
Applying 
Definition \ref{def:pi} to $\Phi^*$ 
we find
$p^*_1=\lambda -a^*_0$.
Combining these facts we get the result.
\hfill $\Box $ \\

\begin{lemma}
\label{lem:abcthsrec}
Let $\Phi$ denote the Leonard system
from 
Definition 
\ref{eq:ourstartingpt} and assume $d\geq 1$.
Then
\begin{eqnarray*}
b_i\theta^*_{i+1} 
+a_i\theta^*_i +
c_i\theta^*_{i-1}
=\theta_1\theta^*_i+a^*_0(\theta_0-\theta_1)
\qquad (0 \leq i\leq d),
\label{eq:abcreq}
\end{eqnarray*}
 where $\theta^*_{-1}, \theta^*_{d+1}$
denote indeterminates.
\end{lemma}
\noindent {\it Proof:}
Set $j=1$ in
(\ref{eq:3termu}).
Evaluate the result using
Lemma
\ref{def:bici}(ii) and 
(\ref{eq:thsform}).
\hfill $\Box $ \\

\section{A bilinear form}

\medskip
\noindent
In this section we associate with each Leonard pair a certain
bilinear form. To prepare for this   
we recall a few concepts from linear
algebra.

\medskip
\noindent 
Let $V$ denote a finite dimensional
vector space over $\K$.
By a {\it bilinear form on $V$} we mean
a  map $\langle \,,\,\rangle : V \times V \rightarrow \K$
that satisfies the following four conditions for all
$u,v,w \in V$ and for all $\alpha \in \K$:
(i) $\langle u+v,w \rangle = 
 \langle u,w \rangle 
+
\langle v,w \rangle$;
(ii)
$
\langle \alpha u,v \rangle 
=
\alpha \langle u,v \rangle 
$;
(iii)
 $\langle u,v+w \rangle = 
 \langle u,v \rangle 
+
\langle u,w \rangle$;
(iv)
$
\langle u, \alpha v \rangle 
=
\alpha \langle u,v \rangle 
$.
We observe 
that a scalar multiple of 
a bilinear form on $V$ is a bilinear form on $V$.
Let
$
\langle\, ,\, \rangle 
$
denote a bilinear form on $V$.
This form is said to be {\it symmetric}
whenever
$
\langle u,v \rangle 
= 
\langle v,u \rangle 
$
for all $u, v \in V$.
Let 
$
\langle\, ,\, \rangle 
$
denote a bilinear form on $V$.
Then the following
are equivalent: (i) there exists a nonzero $u \in V$ such
that $
\langle u,v \rangle 
= 0 
$
for all $v \in V$;
 (ii) there exists a nonzero $v \in V$ such
that $
\langle u,v \rangle 
= 0 
$
for all $u \in V$.
The form 
$\langle \,,\,\rangle $
is 
said to be {\it degenerate }
whenever (i), (ii) hold and {\it nondegenerate}
otherwise.
Let $\gamma :{\mathcal A} \rightarrow {\mathcal A}$
denote an antiautomorphism
and let $V$ denote an irreducible $\mathcal A$-module.
Then there exists a nonzero bilinear form 
$\langle \,,\,\rangle $ on $V$ such that
$\langle Xu,v\rangle 
=
\langle u,X^\gamma v\rangle 
$
for all $u,v \in V$ and for all $X \in \mathcal A$.
The form is unique
up to multiplication by a nonzero scalar in $\K$.
The form in nondegenerate.
We refer to this form as the {\it bilinear form on $V$ associated
with $\gamma $}. This form is not symmetric in general.

\medskip
\noindent We now return our attention to Leonard pairs.

\begin{definition}
\label{def:ip}
\rm
Let
$\Phi=(A;A^*; \lbrace E_i\rbrace_{i=0}^d; $ $ 
\lbrace E^*_i\rbrace_{i=0}^d)$
denote a Leonard system in
$\mathcal A$.
Let
$\dagger :{\mathcal A}\rightarrow {\mathcal A} $ 
denote the corresponding antiautomorphism 
from Definition
\ref{def:dag}.
Let $V$ denote an 
irreducible $\cal A$-module.
For the rest of this paper  we let  
$\langle \,,\,\rangle $ denote the bilinear
form on $V$ 
associated with $\dagger$. 
We abbreviate $\Vert u \Vert^2 = \langle u,u\rangle $
for all $u \in V$.
By the construction,
 for $X \in {\mathcal A}$  we have
\begin{eqnarray}
\langle Xu, v \rangle = \langle u, X^\dagger v \rangle
\qquad (\forall u \in V, \forall v \in V).
\label{eq:daggerstart}
\end{eqnarray}
\end{definition}

\noindent We make an observation.

\begin{lemma}
With reference to Definition \ref{def:ip},
let $\mathcal D$ (resp. $\mathcal D^*$) denote
the subalgebra of $\mathcal A$ generated  by $A$ (resp. $A^*$.)
Then for $X \in {\mathcal D} \cup {\mathcal D}^*$ we have
\begin{eqnarray}
\label{eq:bilinv}
\langle Xu, v \rangle = \langle u, Xv \rangle
\qquad (\forall u \in V, \forall v \in V).
\end{eqnarray}
\end{lemma}
\noindent {\it Proof:}
Combine  (\ref{eq:daggerstart})
and
Lemma \ref{cor:eistab}.
\hfill $\Box $ \\

\medskip
\noindent 
With reference to Definition
\ref{def:ip}, our next goal is to show
$\langle \,,\,\rangle $ is symmetric.
We will use the following lemma.


\begin{theorem}
\label{thm:orthog}
With reference to Definition
\ref{def:ip},
let $u$ denote a nonzero vector in $E_0V$ and
recall 
$E^*_0u, E^*_1u,\ldots, E^*_du $
is a $\Phi$-standard basis for $V$.
We have
\begin{eqnarray}
\langle E^*_iu, E^*_ju \rangle = \delta_{ij}k_i \nu^{-1}   
\Vert u \Vert^2  
\qquad \qquad 
(0 \leq i,j \leq d),
\label{eq:vert}
\end{eqnarray}
where the $k_i$ are from
Definition \ref{def:ki1}
and $\nu$ is from
Definition \ref{def:n}.
\end{theorem}
\noindent {\it Proof:}
By 
(\ref{eq:bilinv}) and since $E_0u=u$ we find
$\langle E^*_iu, E^*_ju\rangle
=
\langle u, E_0E^*_iE^*_jE_0u\rangle$.
Using
Lemma \ref{lem:mid}(ii) and
(\ref{eq:ki1}) we find
$\langle u, E_0E^*_iE^*_jE_0u\rangle=
\delta_{ij}k_i \nu^{-1}   
\Vert u \Vert^2 $. 
\hfill $\Box $ \\

\begin{corollary}
\label{cor:sym}
With reference to
Definition
\ref{def:ip}, the bilinear form 
 $\langle \,, \,\rangle $ is symmetric.
\end{corollary}
\noindent {\it Proof:}
Let $u$ denote a nonzero vector in $E_0V$ and abbreviate
$v_i=E^*_iu$ for $0 \leq i \leq d$.
From Theorem
\ref{thm:orthog} we find
$\langle v_i,v_j\rangle= \langle v_j,v_i\rangle$
for $0 \leq i,j\leq d$. The result follows since
$v_0, v_1,\ldots, v_d$ is a basis for $V$.
\hfill $\Box $ \\

\noindent We have a comment.

\begin{lemma}
\label{lem:anglebasic}
With reference to Definition
\ref{def:ip}, let $u$ denote a nonzero vector in 
$E_0V$ and let  $v$ denote a nonzero vector in $E^*_0V$. Then 
the following (i)--(iv) hold.
\begin{enumerate}
\item Each of
 $\Vert u \Vert^2, 
\Vert v \Vert^2,
\langle u,v\rangle $
is nonzero.
\item
$E^*_0u = 
\langle u,v\rangle 
\Vert v\Vert^{-2} 
 v.$
\item
$
E_0v = 
\langle u,v\rangle 
\Vert u\Vert^{-2} u.$ 
\item $\nu 
\langle u,v\rangle^2 =
\Vert u\Vert^2
\Vert v\Vert^2
$.
\end{enumerate}
\end{lemma}
\noindent {\it Proof:}
(i)
Observe $\Vert u \Vert^2\not=0$ by
Theorem \ref{thm:orthog}
and since $\langle \,,\,\rangle $ is not 0.
Similarly
$\Vert v \Vert^2\not=0$.
To see that $\langle u,v\rangle \not=0$, observe that
$v$ is a basis for $E^*_0V$ so there exists $\alpha \in \K$ 
such that
$E^*_0u = \alpha v$.
Recall $E^*_0u\not=0$ by
Lemma
\ref{lem:stbasis} so $\alpha\not=0$.
Using (\ref{eq:bilinv}) and $E^*_0v=v$ we routinely find
$\langle u,v \rangle =
\alpha \Vert v\Vert^2$ and it follows
$\langle u,v \rangle \not=0$. 
\\
(ii) In the proof of part (i) we 
found $E^*_0u=\alpha v$ where $\langle u,v\rangle =\alpha \Vert v\Vert^2$.
The result follows.
\\
(iii) Similar to the proof of (ii) above.
\\
(iv)
Using $u=E_0u$ and $\nu E_0E^*_0E_0=E_0$ we find  
$\nu^{-1} u = E_0E^*_0u$.
To finish the proof,
evaluate $E_0E^*_0u$ using 
(ii) above and then (iii) above.
\hfill $\Box $ \\

%
%

\begin{theorem}
\label{lem:predu}
With reference to Definition
\ref{def:ip}, let $u$ denote a nonzero vector in 
$E_0V$ and let  $v$ denote a nonzero vector in $E^*_0V$. Then 
\begin{equation}
\langle E^*_iu, E_jv\rangle = 
\nu^{-1}k_i k^*_j u_i(\theta_j)  \langle u,v \rangle 
\qquad  \qquad 
 (0 \leq i,j\leq d).
\label{eq:predu}
\end{equation}
\end{theorem}
\noindent {\it Proof:}
Using Theorem
\ref{def:vi3} we find
\begin{eqnarray}
\langle E^*_iu, E_jv\rangle &=&
\langle v_i(A)E^*_0u, E_jv\rangle 
\nonumber
\\
&=&
\langle E^*_0u,v_i(A)E_jv\rangle 
\nonumber
\\
&=&
v_i(\theta_j)\langle E^*_0u,E_jv\rangle 
\nonumber
\\
&=&
v_i(\theta_j)\langle E^*_0u,v^*_j(A^*)E_0v\rangle 
\nonumber
\\
&=&
v_i(\theta_j)\langle v^*_j(A^*)E^*_0u,E_0v\rangle 
\nonumber
\\
&=&
v_i(\theta_j)v^*_j(\theta^*_0)\langle E^*_0u,E_0v\rangle.
\label{eq:almostdone}
\end{eqnarray}
Using Lemma
\ref{lem:anglebasic}(ii)--(iv) we find
$\langle E^*_0u,E_0v\rangle =\nu^{-1}\langle u,v\rangle$.
Observe
$v_i(\theta_j)=u_i(\theta_j)k_i$
by 
(\ref{eq:uvk}). Applying
Lemma \ref{def:vi1a} to $\Phi^*$ we find
$v^*_j(\theta^*_0)=k^*_j$. 
Evaluating (\ref{eq:almostdone}) using
these comments we obtain
(\ref{eq:predu}).
\hfill $\Box $ \\

\begin{remark} \rm
Using Theorem
\ref{lem:predu} and the symmetry of $\langle \,,\,\rangle $
we get an alternate proof of
Theorem \ref{lem:uidual}.
\end{remark}

\begin{theorem}
\label{thm:Ptrans}
With reference to Definition
\ref{def:ip}, let $u$ denote a nonzero vector in 
$E_0V$ and let  $v$ denote a nonzero vector in $E^*_0V$.
Then for $0 \leq i \leq d$,  both 
\begin{eqnarray}
E^*_i u &=& 
\frac
{
\langle u,v 
\rangle
}
{
\Vert v 
\Vert^2
}
\sum_{j=0}^d v_i(\theta_j) E_jv,
\label{eq:trans1}
\\
E_i v &=& 
\frac
{
\langle u,v 
\rangle
}
{
\Vert u
\Vert^2
}
\sum_{j=0}^d v^*_i(\theta^*_j) E^*_ju.
\label{eq:trans2}
\end{eqnarray}
\end{theorem}
\noindent {\it Proof:}
We first show 
(\ref{eq:trans1}).
To do this we show
each side of 
(\ref{eq:trans1}) is equal to
$v_i(A)E^*_0u$.
By
Theorem \ref{def:vi3} we find
$v_i(A)E^*_0u$ 
is equal to the left-hand side of 
(\ref{eq:trans1}).
To see that 
$v_i(A)E^*_0u$ is equal to the right-hand side
of 
(\ref{eq:trans1}), multiply
$v_i(A)E^*_0u$ on the left by the identity $I$,
expand using
$I=\sum_{j=0}^d E_j$, and simplify the result using
$E_jA=\theta_jE_j$ $(0 \leq j\leq d)$ and
Lemma \ref{lem:anglebasic}(ii).
We have now proved 
(\ref{eq:trans1}).
Applying 
(\ref{eq:trans1}) to $\Phi^*$ we obtain
(\ref{eq:trans2}).
\hfill $\Box $ \\

\begin{definition}
\label{def:pmat}
\rm
Let $\Phi$ denote the Leonard system
from
Definition \ref{eq:ourstartingpt}.
We define a matrix $P \in
\mbox{Mat}_{d+1}(\K)$  as follows.
For $0 \leq i,j\leq d$ the  entry
$P_{ij}=v_j(\theta_i)$, where
$\theta_i$ is from
Definition \ref{def:evseq} and
 $v_j$ is from
Definition
\ref{def:vi1}.
\end{definition}

\begin{theorem}
\label{def:pmatpmats}
Let $\Phi$ denote the Leonard system
from
Definition \ref{eq:ourstartingpt}.
Let the matrix $P$ be as in
Definition
\ref{def:pmat} and recall $P^*$ is the corresponding
matrix for $\Phi^*$.
Then 
$P^*P=\nu I$,
where $\nu$ is from
Definition \ref{def:n}. 
\end{theorem}
\noindent {\it Proof:}
Compare
(\ref{eq:trans1}), (\ref{eq:trans2}) and use
Lemma
\ref{lem:anglebasic}(iv).
\hfill $\Box $ \\

\begin{theorem}
\label{lem:flatvssharp}
Let $\Phi$ denote the Leonard system
from
Definition \ref{eq:ourstartingpt} and
let the matrix $P$ be as in
Definition
\ref{def:pmat}.
Let 
the  map $\flat : {\mathcal A}\rightarrow 
\mbox{Mat}_{d+1}(\K)$ be as in
Definition
\ref{def:flatcon}
and let 
 $\sharp : {\mathcal A}\rightarrow 
\mbox{Mat}_{d+1}(\K)$ denote the corresponding
map for $\Phi^*$.
Then for all $X \in {\mathcal A}$ we have
\begin{eqnarray}
X^\sharp P = P X^\flat.
\end{eqnarray}
\end{theorem}
\noindent {\it Proof:}
Let $V$ denote an irreducible $\mathcal A$-module. 
Let $u$ denote a nonzero vector in $E_0V$
and recall $E^*_0u, E^*_1u, \ldots, E^*_du$ is
a $\Phi$-standard basis for $V$.
By Definition
\ref{def:flatcon}, $X^\flat$ is the matrix in
$\mbox{Mat}_{d+1}(\K)$
that represents 
$X$ with respect to
$E^*_0u, E^*_1u, \ldots, E^*_du$.
Similarly for a nonzero $v \in E^*_0V$,
$X^\sharp$ is the matrix in
$\mbox{Mat}_{d+1}(\K)$
that represents $X$ with respect to
$E_0v, E_1v, \ldots, E_dv$.
In view of (\ref{eq:trans1}),
the 
transition matrix from
$E_0v, E_1v, \ldots, E_dv$ 
to
$E^*_0u, E^*_1u, \ldots, E^*_du$ is
a scalar multiple of $P$.
The result follows 
from these comments and elementary linear 
algebra.
\hfill $\Box $ \\

\section{The orthogonality relations}

\noindent 
In this section we show that
each of the polynomial sequences $p_i, u_i, v_i$
satisfy an orthogonality relation.
We begin with the $v_i$.

\begin{theorem}
\label{thm:viorth}
Let $\Phi$ denote the Leonard system from Definition  
\ref{eq:ourstartingpt} and let
 the polynomials $v_i$ be as in
Definition
\ref{def:vi1}.
Then both 
\begin{eqnarray}
\sum_{r=0}^d v_i(\theta_r)v_j(\theta_r)k^*_r &=&
\delta_{ij} \nu k_i \qquad \qquad (0 \leq i,j\leq d),
\label{eq:vio1}
\\
\sum_{i=0}^d v_i(\theta_r)v_i(\theta_s)k^{-1}_i &=&
\delta_{rs} \nu k_r^{*-1}   \qquad \qquad (0 \leq r,s\leq d).
\label{eq:vio2}
\end{eqnarray}
\end{theorem}
\noindent {\it Proof:}
We refer to Theorem
\ref{def:pmatpmats}.
To obtain 
(\ref{eq:vio1}) compute the $ij$ entry
in
$P^*P=\nu I$
 using matrix multiplication
and evaluate the result
using
Theorem \ref{lem:vidual}.
To obtain 
(\ref{eq:vio2}) compute the $ij$ entry
of $PP^*=\nu I$ using matrix multiplication 
and evaluate the result
using
Theorem \ref{lem:vidual}.
\hfill $\Box $ \\

\noindent 
We now 
turn to the polynomials $u_i$.

\begin{theorem}
\label{thm:uiorth}
Let $\Phi$ denote the Leonard system from 
Definition \ref{eq:ourstartingpt} and
let the polynomials $u_i$ be as in
Definition
\ref{def:ui1}.
Then both 
\begin{eqnarray*}
\sum_{r=0}^d u_i(\theta_r)u_j(\theta_r)k^*_r &=&
\delta_{ij}\nu k^{-1}_i \qquad \qquad (0 \leq i,j\leq d),
\\
\sum_{i=0}^d u_i(\theta_r)u_i(\theta_s)k_i &=&
\delta_{rs} \nu k_r^{*-1}   \qquad \qquad (0 \leq r,s\leq d).
\end{eqnarray*}
\end{theorem}
\noindent {\it Proof:}
Evaluate each of 
(\ref{eq:vio1}), 
(\ref{eq:vio2}) using
Lemma
\ref{def:ui4}.
\hfill $\Box $ \\

\noindent 
We now turn to the polynomials
$p_i$.

\begin{theorem}
\label{thm:piorth}
Let $\Phi$ denote the Leonard system from Definition  
\ref{eq:ourstartingpt}
and let the polynomials $p_i$ be as in
Definition
\ref{def:pi}.
Then both
\begin{eqnarray*}
 \sum_{r=0}^d p_i(\theta_r)p_j(\theta_r)m_r &=& \delta_{ij}
x_1x_2\cdots x_i \qquad \qquad (0 \leq i,j\leq d),
\\
\sum_{i=0}^d {{ p_i(\theta_r) p_i(\theta_s)}\over {x_1x_2\cdots x_i}}
&=& \delta_{rs}m_r^{-1} \qquad \qquad (0 \leq r,s\leq d).
\end{eqnarray*}
\end{theorem}
\noindent {\it Proof:}
Applying Definition
\ref{def:ki1} to $\Phi^*$ we find
$k^*_r=m_r\nu$ for $0 \leq r \leq d$.
Evaluate 
each of
(\ref{eq:vio1}), 
(\ref{eq:vio2}) using
this and 
Definition \ref{def:vi1},
Lemma \ref{def:bici}(i),
(\ref{eq:kibc}).
\hfill $\Box $ \\

\section{Everything in terms of the parameter array}

\noindent 
In this section we express  all the polynomials and
scalars that came up so far in the paper, in terms
of a short list of parameters called the {\it parameter array}. 
 The parameter array of a Leonard system
 consists of
its eigenvalue sequence,
its dual eigenvalue sequence,
and two additional sequences called the {\it first split sequence}
and the {\it  second split sequence}.
The first split sequence  is defined as follows.
Let $\Phi$ denote the Leonard system
from Definition
\ref{eq:ourstartingpt}.
We showed in
\cite[Theorem 3.2]{LS99}
that there exists nonzero scalars $\varphi_1, \varphi_2, \ldots, \varphi_d$
in $\K$ and there exists an isomorphism
of $\K$-algebras
$\natural : {\mathcal A}\rightarrow 
\hbox{Mat}_{d+1}(\K)$
such that
\begin{equation}
A^\natural =\left(
\begin{array}{c c c c c c}
\theta_0 & & & & & {\bf 0} \\
1 & \theta_1 &  & & & \\
& 1 & \theta_2 &  & & \\
& & \cdot & \cdot &  &  \\
& & & \cdot & \cdot &  \\
{\bf 0}& & & & 1 & \theta_d
\end{array}
\right),
\qquad  \quad 
A^{*\natural} = 
\left(
\begin{array}{c c c c c c}
\theta^*_0 &\varphi_1 & & & & {\bf 0} \\
 & \theta^*_1 & \varphi_2 & & & \\
&  & \theta^*_2 & \cdot & & \\
& &  & \cdot & \cdot &  \\
& & &  & \cdot & \varphi_d \\
{\bf 0}& & & &  & \theta^*_d
\end{array}
\right) ,
\label{eq:matrepaastar}
\end{equation}
where the $\theta_i, \theta^*_i$ are from
Definition
\ref{def:evseq}.
The sequence
$\natural, \varphi_1, \varphi_2,\ldots, \varphi_d$
is uniquely determined by $\Phi$.
We call the sequence 
$ \varphi_1, \varphi_2,\ldots, \varphi_d$
the {\it first split sequence} of $\Phi$.
We let 
$ \phi_1, \phi_2,\ldots, \phi_d$
denote the first split sequence
of $\Phi^\Downarrow$ and call
this the 
 {\it second split sequence} of $\Phi$.
For notational convenience we define
$\varphi_0=0$, 
$\varphi_{d+1}=0$, 
$\phi_0=0$, 
$\phi_{d+1}=0$.

\begin{definition}
\label{def:pararray}
\rm
Let $\Phi$ denote the Leonard system
from 
Definition
\ref{eq:ourstartingpt}.
By the {\it parameter array} of $\Phi$ we mean
the sequence
$(\theta_i,\theta^*_i,i=0..d;\varphi_j,\phi_j,j=1..d)$,
where  
 $\theta_0, \theta_1, \ldots, \theta_d$ 
  (resp. $\theta^*_0, \theta^*_1, \ldots, \theta^*_d$)
  is the eigenvalue sequence (resp. 
  dual eigenvalue sequence) of $\Phi$
and   
$ \varphi_1, \varphi_2,\ldots, \varphi_d$
(resp.
$ \phi_1, \phi_2,\ldots, \phi_d$)
is 
the first split sequence
(resp. second split sequence) of
 $\Phi$.
\end{definition}

\noindent We remark that two
Leonard systems over $\K$ are isomorphic if and only if they have the
same parameter array
\cite[Theorem 1.9]{LS99}.

\medskip
\noindent 
The following result shows that 
the parameter array
 behaves nicely with respect to the $D_4$ action given in
 Section 4.

\begin{theorem} \cite[Theorem 1.11]{LS99}
\label{thm:phimod}
Let $\Phi$ denote a Leonard system
with parameter array
$(\theta_i,\theta^*_i,i=0..d; $ $\varphi_j,\phi_j,j=1..d)$.
 Then  (i)--(iii) hold below.
\begin{enumerate}
\item  The parameter array of $\Phi^*$
is
$(\theta^*_i,\theta_i,i=0..d;\varphi_j,\phi_{d-j+1},j=1..d)$.
\item 
 The parameter array of
 $\Phi^{\downarrow}$
is
$(\theta_i,\theta^*_{d-i},i=0..d;\phi_{d-j+1},\varphi_{d-j+1},j=1..d)$.
\item  The parameter array of $\Phi^\Downarrow$
is
$(\theta_{d-i},\theta^*_i,i=0..d;\phi_j,\varphi_j,j=1..d)$.
\end{enumerate}
\end{theorem}

\noindent For the rest of this paper we will 
use the following notation.

\begin{definition}
\label{def:tau}
\rm
Suppose we are given an integer
$d\geq 0$ and two sequences of scalars
\begin{eqnarray*}
\theta_0, \theta_1, \ldots, \theta_d; \quad
\theta^*_0, \theta^*_1, \ldots, \theta^*_d
\end{eqnarray*}
taken from $\K$.
Then for $0 \leq i \leq d+1$ we let 
$\tau_i$, 
$\tau^*_i$, 
$\eta_i$,
$\eta^*_i$
denote the following polynomials in
$\K \lbrack \lambda \rbrack $.
\begin{eqnarray}
&&\tau_i = \prod_{h=0}^{i-1} (\lambda - \theta_h),
\qquad \qquad \;\;
\tau^*_i = \prod_{h=0}^{i-1} (\lambda - \theta^*_h),
\label{eq:tau}
\label{eq:taus}
\\
&&\eta_i = \prod_{h=0}^{i-1}(\lambda - \theta_{d-h}),
\qquad \qquad
\eta^*_i = \prod_{h=0}^{i-1}(\lambda - \theta^*_{d-h}).
\label{eq:eta}
\label{eq:etas}
\end{eqnarray}
\noindent 
We observe that each of
$\tau_i$, 
$\tau^*_i$, 
$\eta_i$,
$\eta^*_i$
is monic with degree $i$.
\end{definition}

\begin{theorem}
\label{thm:uform}
Let $\Phi$ 
denote the Leonard system from Definition
\ref{eq:ourstartingpt}
and 
let 
$(\theta_i, \theta^*_i, i=0..d;
\varphi_j,\phi_j,j=1..d)$ denote the corresponding
parameter array.
Let the polynomials $u_i$ be as in
Definition \ref{def:ui1}.
Then
\begin{eqnarray}
\label{eq:uifinal}
u_i=\sum_{h=0}^i {{\tau^*_h(\theta^*_i)}\over 
{\varphi_1\varphi_2\cdots \varphi_h}} \tau_h 
\qquad \qquad 
(0 \leq i \leq d).
\end{eqnarray}
We are using the notation
 (\ref{eq:taus}).
\end{theorem}
\noindent {\it Proof:}
Let the integer $i$ be given. 
The polynomial $u_i$ has degree $i$ so
 there  exists 
 scalars $\alpha_0, \alpha_1, \ldots, \alpha_i$ 
in $\K$ such that
\begin{eqnarray}
\label{eq:ualphaform}
u_i = \sum_{h=0}^i \alpha_h \tau_h.
\end{eqnarray}
We show
\begin{eqnarray}
\alpha_h = \frac{\tau^*_h(\theta^*_i)}{\varphi_1\varphi_2 \cdots \varphi_h}
\qquad \qquad (0 \leq h \leq i).
\label{eq:alphagoal}
\end{eqnarray}
In order to do this we show
$\alpha_0=1$ and 
$\alpha_{h+1}\varphi_{h+1}=\alpha_h(\theta^*_i-\theta^*_h)$
for $0 \leq h \leq i-1$.
We now show $\alpha_0=1$.
We evaluate (\ref{eq:ualphaform})
at $\lambda=\theta_0$
and find $u_i(\theta_0)=\sum_{h=0}^i \alpha_h \tau_h(\theta_0)$.
Recall $u_i(\theta_0)=1$ by
(\ref{eq:unorm}).
Using 
(\ref{eq:tau})
we find $\tau_h(\theta_0)=1$  for $h=0$ and
$\tau_h(\theta_0)=0$  for
$1 \leq h \leq i$.
From these comments  
we find $\alpha_0=1$.
We  now show
$\alpha_{h+1} \varphi_{h+1} = 
\alpha_h(\theta^*_i-\theta^*_h)$
for $0 \leq h \leq i-1$.
Let $V$ denote an irreducible $\mathcal A$-module.
From 
(\ref{eq:matrepaastar}) there exists a basis
 $e_0, e_1, \ldots, e_d$ for $V$
that satisfies
$(A-\theta_jI)e_j = e_{j+1}$ $(0 \leq j \leq d-1)$,
$(A-\theta_dI)e_d =0$
and $(A^*-\theta^*_jI)e_j = \varphi_je_{j-1}$ 
$(1 \leq j \leq d)$,
$(A^*-\theta^*_0I)e_0 =0$.
From the action of $A$ on $e_0, e_1, \ldots, e_d$
we find $e_j = \tau_j(A)e_0$ for $0 \leq j\leq d$.
Observe $A^*e_0=\theta^*_0e_0$ so $e_0 \in E^*_0V$.
Combining 
Theorem \ref{thm:psend}
and 
(\ref{uicl}) we find
$u_i(A)E^*_0V=E^*_iV$.
By this and since
 $e_0 \in E^*_0V$ 
we find
$u_i(A)e_0 \in E^*_iV$. 
Apparently 
$u_i(A)e_0 $ is an eigenvector for $A^*$ with
eigenvalue $\theta^*_i$.
We may now argue
\begin{eqnarray*}
0  &=& (A^*-\theta^*_iI)u_i(A)e_0
\\
  &=& (A^*-\theta^*_iI)\sum_{h=0}^i \alpha_h \tau_h(A)e_0
\\
  &=& (A^*-\theta^*_iI)\sum_{h=0}^i \alpha_h e_h 
\\
  &=& \sum_{h=0}^{i-1} e_h 
  (\alpha_{h+1} \varphi_{h+1}-\alpha_{h}(\theta^*_i-\theta^*_h)).
\end{eqnarray*}
By this and since $e_0, e_1, \ldots, e_d$ are linearly independent we find
  $\alpha_{h+1} \varphi_{h+1}=\alpha_h(\theta^*_i-\theta^*_h)$ for
  $0 \leq h\leq i-1$.
Line
(\ref{eq:alphagoal}) follows and the theorem is proved.
\hfill $\Box $ \\

\begin{lemma}
\label{lem:p0}
Let $\Phi$ denote the Leonard system
from 
Definition 
\ref{eq:ourstartingpt}
and let
$(\theta_i, \theta^*_i, i=0..d;
\varphi_j,\phi_j,j=1..d)$ denote the corresponding
parameter array.
Let the polynomials 
$p_i$ be as in
Definition
\ref{def:pi}.
With reference to
Definition \ref{def:tau} we have
\begin{eqnarray}
\label{eq:pith0first}
p_i(\theta_0)=  \frac{\varphi_1\varphi_2\cdots \varphi_i}{\tau^*_i(\theta^*_i)}
\qquad \qquad (0 \leq i \leq d).
\end{eqnarray}
\end{lemma}
\noindent {\it Proof:}
In equation
(\ref{eq:uifinal}), each side is a polynomial
of degree $i$ in $\lambda $.
For the polynomial on the left 
in 
(\ref{eq:uifinal})
the coefficient of $\lambda^i$
is $p_i(\theta_0)^{-1}$
by 
(\ref{uicl}) and since  $p_i$ is monic.
For the polynomial
on the right
in 
(\ref{eq:uifinal})
the coefficient of
 $\lambda^i$ is
$\tau^*_i(\theta^*_i)(\varphi_1\varphi_2\cdots \varphi_i)^{-1}$.
Comparing these coefficients we obtain
the result.
\hfill $\Box $ \\

\begin{theorem}
\label{thm:pform}
Let $\Phi$ 
denote the Leonard system from Definition
\ref{eq:ourstartingpt}
and 
let 
$(\theta_i, \theta^*_i, i=0..d;
\varphi_j,\phi_j,j=1..d)$ denote the corresponding
parameter array.
Let the polynomials 
$p_i$ be as in
Definition \ref{def:pi}.
Then with reference to
Definition \ref{def:tau},
\begin{eqnarray*}
p_i=\sum_{h=0}^i \frac{\varphi_1\varphi_2\cdots \varphi_i}{
\varphi_1\varphi_2\cdots \varphi_h}
\frac{\tau^*_h(\theta^*_i)} 
{\tau^*_i(\theta^*_i)} \tau_h
\qquad \qquad 
(0 \leq i \leq d).
\end{eqnarray*}
\end{theorem}
\noindent {\it Proof:}
Observe $p_i = p_i(\theta_0)u_i$ by
(\ref{uicl}). In this equation we evaluate
$p_i(\theta_0)$ using
(\ref{eq:pith0first}) and
we evaluate $u_i$ using
(\ref{eq:uifinal}). The result follows.
\hfill $\Box $ \\

\begin{theorem}
\label{thm:bcform}
Let $\Phi$ 
denote the Leonard system from Definition
\ref{eq:ourstartingpt}
and 
let 
$(\theta_i, \theta^*_i, i=0..d;
\varphi_j,\phi_j,j=1..d)$ denote the corresponding
parameter array.
Let the scalars  $b_i, c_i$ be as in
Definition
\ref{def:sharpmp}.
Then with reference to
Definition \ref{def:tau}
the following (i), (ii) hold.
\begin{enumerate}
\item
$\displaystyle
{b_i = \varphi_{i+1} \frac{\tau^*_i(\theta^*_i)}{
\tau^*_{i+1}(\theta^*_{i+1})} \qquad \qquad 
(0 \leq i \leq d-1)
}$.
\\
\item
$
\displaystyle{
c_i = \phi_i \frac{\eta^*_{d-i}(\theta^*_i)}{\eta^*_{d-i+1}(\theta^*_{i-1})}
\qquad \qquad (1 \leq i \leq d)}$.
\end{enumerate}
\end{theorem}
\noindent {\it Proof:}
(i) Evaluate
(\ref{eqbi}) using
Lemma
\ref{lem:p0}.
\\
(ii) Comparing
the formulae for $b_i, c_i$ given in
Theorem \ref{def:bici4}
we find, with reference to
Definition
\ref{def:notconv},
that $c_i=b^\downarrow_{d-i}$.
Applying
part (i) above to $\Phi^\downarrow$
and using
Theorem \ref{thm:phimod}(ii) we routinely
obtain the result.
\hfill $\Box $ \\

\noindent Let $\Phi$ denote the Leonard system
from 
 Definition
\ref{eq:ourstartingpt}
and let the scalars $a_i$ be as in
Definition
\ref{def:aixi}.
We  mention two formulae that give  
$a_i$ in terms of the parameter array of $\Phi$.
The first formula
is obtained using
Lemma
\ref{def:bici}(ii) and
Theorem \ref{thm:bcform}.
The second formula is given in the following theorem.
This theorem was proven in
\cite[Lemma 5.1]{LS99}; however we give an alternate proof
that we find illuminating.

\begin{theorem}
\cite[Lemma 5.1]{LS99}
\label{thm:ai}
Let $\Phi$ 
denote the Leonard system from Definition
\ref{eq:ourstartingpt}
and 
let 
$(\theta_i, \theta^*_i, i=0..d;
\varphi_j,\phi_j,j=1..d)$ denote the corresponding
parameter array.
Let the scalars $a_i$ be as in 
Definition
\ref{def:aixi}. Then
\begin{eqnarray}
a_i = \theta_i + \frac{\varphi_i}{\theta^*_i-\theta^*_{i-1}}
+ \frac{\varphi_{i+1}}{\theta^*_i-\theta^*_{i+1}} \qquad (0 \leq i \leq d),
\label{eq:aiform}
\end{eqnarray}
where we recall $\varphi_0=0$, $\varphi_{d+1}=0$, and where
$\theta^*_{-1}, \theta^*_{d+1}$ denote indeterminates.
\end{theorem}
\noindent {\it Proof:}
Let the polynomials $p_0, p_1, \ldots, p_{d+1}$ be as in
Definition
\ref{def:pi} and recall these polynomials are monic.
Let $i$ be given and consider the polynomial
\begin{eqnarray}
\lambda p_i-p_{i+1}.
\label{eq:picomb}
\end{eqnarray}
From (\ref{eq:prec}) 
we find
the polynomial (\ref{eq:picomb}) is equal to 
$a_ip_i+x_ip_{i-1}$.
Therefore 
the polynomial
(\ref{eq:picomb})
has degree $i$ and leading coefficient $a_i$.
In order to compute this leading coefficient,
in 
(\ref{eq:picomb})
we 
evaluate each of 
$p_i$, $p_{i+1}$
using
Theorem \ref{thm:pseq}(ii)
and 
Theorem
\ref{thm:pform}.
By this method we routinely obtain
(\ref{eq:aiform}).
\hfill $\Box $ \\

\begin{theorem}
\label{thm:xform}
Let $\Phi$ 
denote the Leonard system from Definition
\ref{eq:ourstartingpt}
and 
let 
$(\theta_i, \theta^*_i, i=0..d;
\varphi_j,\phi_j,j=1..d)$ denote the corresponding
parameter array.
Let the scalars $x_i$ be
as in
Definition
\ref{def:aixi}.
Then with reference to
Definition \ref{def:tau},
\begin{eqnarray}
\label{eq:xifinform}
x_i = 
 \varphi_{i}
\phi_i  
 \frac{\tau^*_{i-1}(\theta^*_{i-1})
\eta^*_{d-i}(\theta^*_i)}
{
\tau^*_{i}(\theta^*_{i})
\eta^*_{d-i+1}(\theta^*_{i-1})}
\qquad \qquad (1 \leq i \leq d).
\end{eqnarray}
\end{theorem}
\noindent {\it Proof:}
Use $x_i=b_{i-1}c_i$ and
Theorem
\ref{thm:bcform}.
\hfill $\Box $ \\

\begin{theorem}
\label{thm:nform}
Let $\Phi$ 
denote the Leonard system from Definition
\ref{eq:ourstartingpt}
and 
let 
$(\theta_i, \theta^*_i, i=0..d;
\varphi_j,\phi_j,j=1..d)$ denote the corresponding
parameter array.
Let the scalar $\nu$ be
as in
Definition
\ref{def:n}. 
Then with reference to
Definition \ref{def:tau},
\begin{eqnarray}
\nu = \frac{\eta_d(\theta_0) \eta^*_d(\theta^*_0)
}
{\phi_1 \phi_2 \cdots \phi_d 
}.
\label{eq:nfinalform}
\end{eqnarray}
\end{theorem}
\noindent {\it Proof:}
Evaluate 
(\ref{frame}) using
Theorem \ref{thm:bcform}(ii).
\hfill $\Box $ \\

\begin{theorem}
\label{thm:kiform}
Let $\Phi$ 
denote the Leonard system from Definition
\ref{eq:ourstartingpt}
and 
let 
$(\theta_i, \theta^*_i, i=0..d;
\varphi_j,\phi_j,j=1..d)$ denote the corresponding
parameter array.
Let the scalars $k_i$ be
as in
Definition
\ref{def:ki1}.
Then with reference to
Definition \ref{def:tau},
\begin{eqnarray}
\label{eq:kform2}
k_i = 
\frac{\varphi_1\varphi_2\cdots \varphi_i} 
{\phi_1\phi_2\cdots \phi_i}
\frac{\eta^*_d(\theta^*_0)}{\tau^*_i(\theta^*_i)\eta^*_{d-i}(\theta^*_i)}
\qquad \qquad (0 \leq i \leq d).
\end{eqnarray}
\end{theorem}
\noindent {\it Proof:}
Evaluate 
(\ref{eq:kibc})
using 
Theorem \ref{thm:bcform}.
\hfill $\Box $ \\

\begin{theorem}
\label{thm:miform}
Let $\Phi$ 
denote the Leonard system from Definition
\ref{eq:ourstartingpt}
and 
let 
$(\theta_i, \theta^*_i, i=0..d;
\varphi_j,\phi_j,j=1..d)$ denote the corresponding
parameter array.
Let the scalars $m_i$ be
as in
Definition
\ref{def:mid}.
Then with reference to
Definition \ref{def:tau},
\begin{eqnarray}
\label{eq:mform2}
m_i = 
\frac{\varphi_1\varphi_2\cdots \varphi_i 
\phi_1\phi_2\cdots \phi_{d-i}}
{\eta^*_d(\theta^*_0)\tau_i(\theta_i)\eta_{d-i}(\theta_i)}
\qquad \qquad (0 \leq i \leq d).
\end{eqnarray}
\end{theorem}
\noindent {\it Proof:}
Applying
Definition
\ref{def:ki1} to
$\Phi^*$ we find
$m_i=k^*_i\nu^{-1}$.
We compute $k^*_i$ using
Theorem \ref{thm:kiform}
and Theorem
\ref{thm:phimod}(i).
We compute 
$\nu$ using
Theorem \ref{thm:nform}.
The result follows.
\hfill $\Box $ \\

\section{Some polynomials from the Askey scheme}

Let $\Phi$ denote the Leonard system from
Definition
\ref{eq:ourstartingpt} and let the polynomials
$u_i$ be as in
Definition \ref{def:ui1}. In this section we discuss 
how the $u_i$ fit into the Askey scheme
\cite{KoeSwa}, 
 \cite[p260]{BanIto}. 
Our argument is summarized as follows.
In \cite{TLT:array} 
we displayed 13 families of parameter arrays.
By \cite[Theorem 5.16]{TLT:array}
every parameter array is contained in
at least one of these families. 
In (\ref{eq:uifinal})  the $u_i$ are expressed  as a sum involving
the parameter array of $\Phi$.
In \cite[Examples 5.3-5.15]{TLT:array} we
evaluated this sum 
 for the 13 families 
of parameter arrays.
We found the corresponding $u_i$  
form a class  
consisting of the 
$q$-Racah, $q$-Hahn, dual $q$-Hahn, 
$q$-Krawtchouk,
dual $q$-Krawtchouk,
quantum 
$q$-Krawtchouk,
affine 
$q$-Krawtchouk,
Racah, Hahn, dual Hahn, Krawtchouk,  Bannai/Ito, 
and 
orphan polynomials. 
This class coincides with the terminating branch of the 
Askey scheme.
We remark the Bannai/Ito polynomials can be obtained
from the $q$-Racah polynomials by letting $q$ tend to $-1$ 
\cite[p260]{BanIto}. The orphan polynomials
exist for diameter $d=3$ and $\mbox{Char}(\K)=2$ only 
\cite[Example 5.15]{TLT:array}.
We will not reproduce all the details of our calculations
here; instead we illustrate what is going on with
 some examples.
We will consider two families of parameter  arrays.
For the first family the corresponding 
$u_i$ will turn out to be some Krawtchouk
polynomials. For the  second family the corresponding
$u_i$ will turn out to
be the $q$-Racah polynomials.

\medskip
\noindent Our first example
is associated with the Leonard pair
(\ref{eq:fam1}). 
Let $d$ denote a nonnegative integer and 
consider the following elements of $\K$.
\begin{eqnarray}
&&\theta_i =  d-2i, \qquad \qquad \theta^*_i = d-2i 
\qquad \qquad (0 \leq i \leq d),
\label{eq:thsol}
\\
&&
\varphi_i = -2i(d-i+1), \qquad \qquad \phi_i = 2i(d-i+1) 
\qquad  \qquad (1 \leq i \leq d).
\label{eq:vpsol}
\end{eqnarray}
In order to avoid degenerate situations
we assume the characteristic
of $\K$ is zero or an odd prime greater than $d$.
By 
\cite[Theorem 1.9]{LS99}
%
we find there exists a Leonard system
$\Phi$
over $\K$ that has parameter array
$(\theta_i, \theta^*_i,i=0..d; \varphi_j, \phi_j,j=1..d)$.
Let
the scalars $a_i$ for
  $\Phi$ be as in
  (\ref{eq:aitr}).
Applying Theorem
\ref{thm:ai} to $\Phi$ we find
\begin{eqnarray}
\label{eq:akraw}
a_i=0 \qquad (0 \leq i \leq d).
\end{eqnarray}
Let the scalars $b_i, c_i$ for $\Phi$
be as in
Definition
\ref{def:sharpmp}.
Applying
Theorem 
\ref{thm:bcform} to $\Phi$ we  find
\begin{eqnarray}
\label{eq:bckraw}
b_i = d-i, \qquad c_i = i \qquad \qquad (0 \leq i \leq d).
\end{eqnarray}
Pick any integers $i,j$ $(0 \leq i,j\leq d)$.
Applying 
Theorem
\ref{thm:uform}
to $\Phi$
we find
\begin{eqnarray}
u_i(\theta_j)=\sum_{n=0}^d \frac{(-i)_n (-j)_n 2^n}{(-d)_n n! }, 
\label{eq:2F1expand}
\end{eqnarray}
 where
\begin{eqnarray*}
(a)_n:=a(a+1)(a+2)\cdots (a+n-1) \qquad \qquad n=0,1,2,\ldots
\end{eqnarray*}
Hypergeometric series are defined in \cite[p. 3]{GR}.
From this definition we find
the sum on the right in
(\ref{eq:2F1expand}) is the hypergeometric series
\begin{eqnarray}
{{}_2}F_1\Biggl({{-i, -j}\atop {-d}}\;\Bigg\vert \;2\Biggr).
\label{eq:2F1not}
\end{eqnarray}
A definition of the Krawtchouk polynomials can be found in 
 \cite{AAR} or 
\cite{KoeSwa}. Comparing this definition with 
(\ref{eq:2F1expand}),
(\ref{eq:2F1not})
we find the $u_i$  
are 
Krawtchouk polynomials but not the most general ones. 
Let the scalar $\nu$ for $\Phi$ be as in
Definition \ref{def:n}. 
Applying
Theorem
\ref{thm:nform}
to $\Phi$
 we find 
$\nu=2^d. $
Let the scalars $k_i$ for $\Phi$ be as in
Definition \ref{def:ki1}.
Applying 
Theorem
\ref{thm:kiform}
to $\Phi$
we obtain  
a binomial coefficent
\begin{eqnarray*}
k_i=
\Biggl({{ d }\atop {i}}\Biggr)  \qquad \qquad (0 \leq i\leq d).
\end{eqnarray*}
Let the scalars $m_i$ for $\Phi$ be as in
Definition
\ref{def:mid}.
Applying
Theorem
\ref{thm:miform}
to $\Phi$
we find 
\begin{eqnarray*}
m_i = 
\Biggl({{ d }\atop {i}}\Biggr)2^{-d}  \qquad \qquad (0 \leq i \leq d).
\end{eqnarray*}

\medskip
\noindent
We now give our second example. For this
example the 
 polynomials $u_i$ will turn out to be the
$q$-Racah polynomials. To begin, 
let  $d$ denote a nonnegative integer  and consider the 
following elements in $\K$.
\begin{eqnarray}
\theta_i &=& \theta_0 + h(1-q^i)(1-sq^{i+1})/q^i,
\label{eq:thdefend}
\\
\theta^*_i &=& \theta^*_0 + h^*(1-q^i)(1-s^*q^{i+1})/q^i
\label{eq:thsdefend}
\end{eqnarray}
for $0 \leq i \leq d$, and
\begin{eqnarray}
\varphi_i &=& hh^*q^{1-2i}(1-q^i)(1-q^{i-d-1})(1-r_1q^i)(1-r_2q^i),
\label{eq:varphidefend}
\\
\phi_i &=& hh^*q^{1-2i}(1-q^i)(1-q^{i-d-1})(r_1-s^*q^i)(r_2-s^*q^i)/s^*
\label{eq:phidefend}
\end{eqnarray}
for $1 \leq i \leq d$. We assume 
$q, h, h^*, s, s^*, r_1, r_2$ are nonzero scalars
in the algebraic closure of $\K$, and that $r_1r_2 = s s^*q^{d+1}$.
To avoid degenerate situations
we assume
none of $q^i, r_1q^i, r_2q^i, s^*q^i/r_1, s^*q^i/r_2$ is
equal to 1 for $1 \leq i \leq d$
and neither of
$sq^i, s^*q^i$ is equal to 1 for $2 \leq i \leq 2d$.
By
\cite[Theorem 1.9]{LS99}
there exists a Leonard system $\Phi$ over $\K$ that has 
parameter array
$(\theta_i, \theta^*_i,i=0..d; \varphi_j, \phi_j,j=1..d)$.
Let the scalars $b_i, c_i$ for $\Phi$ be as in
Definition
\ref{def:sharpmp}.
Applying
Theorem
\ref{thm:bcform} to $\Phi$ we find
\begin{eqnarray*}
b_0 &=& \frac{h(1-q^{-d})(1-r_1q)(1-r_2q)}
{1-s^*q^2},
\\
b_i &=& \frac{h(1-q^{i-d})(1-s^*q^{i+1})(1-r_1q^{i+1})(1-r_2q^{i+1})}
{(1-s^*q^{2i+1})(1-s^*q^{2i+2})} \qquad \quad (1 \leq i \leq d-1),
\\
c_i &=& \frac{h(1-q^i)(1-s^*q^{i+d+1})(r_1-s^*q^i)(r_2-s^*q^i)}
{s^*q^d(1-s^*q^{2i})(1-s^*q^{2i+1})} \qquad (1 \leq i \leq d-1),
\\
c_d &=& \frac{h(1-q^d)(r_1-s^*q^d)(r_2-s^*q^d)}
{s^*q^d(1-s^*q^{2d})}.
\end{eqnarray*}
Pick  integers $i,j$ $(0 \leq i,j\leq d)$.
Applying 
Theorem
\ref{thm:uform}
to $\Phi$
we find
\begin{eqnarray}
u_i(\theta_j)=\sum_{n=0}^d \frac{(q^{-i};q)_n (s^*q^{i+1};q)_n 
(q^{-j};q)_n (sq^{j+1};q)_n q^n}
{(r_1q;q)_n(r_2q;q)_n (q^{-d};q)_n(q;q)_n},
\label{eq:uihyper}
\end{eqnarray}
where 
\begin{eqnarray*}
(a;q)_n := (1-a)(1-aq)(1-aq^2)\cdots (1-aq^{n-1})\qquad \qquad n=0,1,2\ldots 
\end{eqnarray*}
Basic hypergeometric series are defined in \cite[p. 4]{GR}.
From that definition we find
the sum on the right in 
(\ref{eq:uihyper}) is the basic hypergeometric series
\begin{eqnarray}
 {}_4\phi_3 \Biggl({{q^{-i}, \;s^*q^{i+1},\;q^{-j},\;sq^{j+1}}\atop
{r_1q,\;\;r_2q,\;\;q^{-d}}}\;\Bigg\vert \; q,\;q\Biggr).
\label{eq:qrac}
\end{eqnarray}
A definition of the  $q$-Racah polynomials can be found in 
\cite{AWil} or 
\cite{KoeSwa}. Comparing this definition with 
(\ref{eq:uihyper}),
(\ref{eq:qrac})
 and
recalling $r_1r_2=s s^*q^{d+1}$,
we find  the $u_i$
are the $q$-Racah polynomials.
Let the scalar $\nu$ for $\Phi$ be as in
Definition \ref{def:n}. 
Applying 
Theorem
\ref{thm:nform}
to $\Phi$
 we find 
\begin{eqnarray*}
\nu = \frac{(sq^2;q)_d (s^*q^2;q)_d}{r^d_1q^d(sq/r_1;q)_d(s^*q/r_1;q)_d}. 
\end{eqnarray*}
Let
the scalars  $k_i$ for $\Phi$ be as in
Definition
\ref{def:ki1}.
Applying
Theorem
\ref{thm:kiform}
to $\Phi$
we obtain  
\begin{eqnarray*}
k_i = \frac{(r_1q;q)_i(r_2q;q)_i(q^{-d};q)_i(s^*q;q)_i(1-s^*q^{2i+1})}
{s^iq^i(q;q)_i(s^*q/r_1;q)_i(s^*q/r_2;q)_i(s^*q^{d+2};q)_i(1-s^*q)} 
\label{eq:qrack}
\qquad (0 \leq i \leq d).
\end{eqnarray*}
Let the scalars $m_i$ for $\Phi$ be as in
Definition \ref{def:mid}.
Applying
Theorem
\ref{thm:miform}
to $\Phi$
we find
\begin{eqnarray*}
m_i = \frac{(r_1q;q)_i(r_2q;q)_i(q^{-d};q)_i(sq;q)_i(1-sq^{2i+1})}
{s^{*i}q^i(q;q)_i(sq/r_1;q)_i(sq/r_2;q)_i(sq^{d+2};q)_i(1-sq)\nu} 
\qquad (0 \leq i \leq d).
\end{eqnarray*}

\section{A characterization of Leonard systems}

\noindent
 In \cite[Appendix A]{LS99} we 
mentioned that the concept of
 a Leonard system
can be viewed as a ``linear algebraic version'' of the polynomial
system which D. Leonard considered 
in \cite{Leon}.
In that appendix we outlined a correspondence that supports this
view but we gave no proof. In this section we provide the proof.

\medskip
\noindent We recall some results from
 earlier in
the paper.
Let $\Phi$ denote the Leonard system from Definition 
\ref{eq:ourstartingpt}.
Let the polynomials
 $p_0, p_1, \ldots, p_{d+1}$
be as in Definition
\ref{def:pi} and recall 
$p^*_0, p^*_1, $ $ \ldots, p^*_{d+1} $
are the corresponding polynomials for $\Phi^*$.
For the purpose of this section, we call
 $p_0, p_1, \ldots, $ $ p_{d+1}$
the {\it monic polynomial sequence} (or {\it MPS}) of $\Phi$.
We call 
$p^*_0, p^*_1,\ldots, p^*_{d+1}$  the  {\it dual MPS} of  $\Phi$.
By Definition 
\ref{def:pi} 
we have
\begin{eqnarray}
&& \qquad p_0=1, \qquad \qquad p^*_0=1, 
\label{eq:rec1introS99}
\\
\lambda p_i &=& p_{i+1} + a_ip_i +x_ip_{i-1} \qquad \qquad (0 \leq i \leq d),
\label{eq:rec2introS99}
\\
\lambda p^*_i &=& p^*_{i+1} + a^*_ip^*_i +x^*_ip^*_{i-1} 
\qquad \qquad (0 \leq i \leq d),
\label{eq:rec3introS99}
\end{eqnarray}
where $\,x_0,\; x^*_0, \; 
p_{-1},\; 
p^*_{-1}\,$ are all zero,  and where 
\begin{eqnarray*} 
 a_i =\mbox{tr}(E^*_iA),\quad && \qquad    \quad   
a^*_i =\mbox{tr}(E_iA^*) \qquad \qquad \quad  (0 \leq i \leq d), \qquad 
\\
x_i =\mbox{tr}(E^*_iAE^*_{i-1}A),&& \qquad     
x^*_i =\mbox{tr}(E_iA^*E_{i-1}A^*) \qquad  \quad   (1 \leq i \leq d). 
\end{eqnarray*}
By 
Lemma \ref{lem:bmat}(iii) we have
\begin{equation}
x_i\not=0, \qquad \qquad 
x^*_i\not=0 \qquad \qquad (1 \leq i \leq d).
\label{eq:rec4introS99}
\end{equation}
Let $\theta_0,\theta_1,\ldots, \theta_d$ 
(resp. $\theta^*_0,\theta^*_1,\ldots,\theta^*_d$)
denote the eigenvalue sequence (resp. dual eigenvalue sequence)
of $\Phi$, and  recall
\begin{eqnarray}
&&\theta_i\not=\theta_j, \qquad \quad 
\theta^*_i\not=\theta^*_j \qquad \quad \hbox{if} \quad i\not=j,  
\qquad \;\;  (0 \leq i,j\leq d).
\label{eq:rec5introS99}
\end{eqnarray}
By 
Theorem
\ref{thm:pseq}(ii) we have
\begin{eqnarray}
&&\qquad p_{d+1}(\theta_i) = 0, \qquad \qquad 
p^*_{d+1}(\theta^*_i) = 0 \qquad  (0 \leq i \leq d). \qquad 
\label{eq:rec6introS99}
\end{eqnarray}
By Theorem 
\ref{def:bici3} we have
\begin{equation} 
p_i(\theta_0)\not=0, \qquad \qquad  
p^*_i(\theta^*_0)\not=0 \qquad \qquad  (0 \leq i \leq d).
\label{eq:rec7introS99}
\end{equation}
By Theorem
\ref{lem:pidual}
we have
\begin{equation}
{{p_i(\theta_j)}\over {
p_i(\theta_0)}}
= 
{{p^*_j(\theta^*_i)}\over {
p^*_j(\theta^*_0)}}\qquad \qquad (0 \leq i,j\leq d). 
\label{eq:rec8introS99}
\end{equation}

\noindent In the following theorem we show the
equations
(\ref{eq:rec1introS99})--(\ref{eq:rec8introS99})
characterize the Leonard systems.

\begin{theorem}
Let $d$ denote a nonnegative integer.
Given 
 polynomials
\begin{eqnarray}
 && p_0, p_1, \ldots, p_{d+1}
\label{eq:172L},
\\
&&p^*_0, p^*_1, \ldots, p^*_{d+1}
\label{eq:172R}
\end{eqnarray}
in 
 $\K\lbrack \lambda \rbrack $ satisfying  
(\ref{eq:rec1introS99})--(\ref{eq:rec4introS99}) 
and given   
 scalars 
\begin{eqnarray}
&&\theta_0, \theta_1, \ldots, \theta_d,
\label{eq:173L}
\\
&&\theta^*_0, \theta^*_1, \ldots, \theta^*_d 
\label{eq:173R}
\end{eqnarray}
in 
 $\K $ 
satisfying 
(\ref{eq:rec5introS99})--(\ref{eq:rec8introS99}),
 there exists  
 a Leonard system  $\Phi$ over $\K$
that has  MPS 
(\ref{eq:172L}),
 dual MPS
(\ref{eq:172R}),
  eigenvalue sequence 
(\ref{eq:173L})
and dual eigenvalue sequence (\ref{eq:173R}).
The system $\Phi $ is unique up to isomorphism of Leonard systems.   
\end{theorem}
\noindent {\it Proof:}
We abbreviate $V={\K}^{d+1}$.
Let $A$ and $A^*$ denote the following matrices in
$\mbox{Mat}_{d+1}(\K)$:
\begin{eqnarray*}
\label{eq:amonmat}
A:=\left(
\begin{array}{ c c c c c c}
a_0 & x_1  &      &      &   &{\bf 0} \\
1 & a_1  &  x_2   &      &   &  \\
  & 1  &  \cdot    & \cdot  &   & \\
  &   & \cdot     & \cdot  & \cdot   & \\
  &   &           &  \cdot & \cdot & x_d \\
{\bf 0} &   &   &   & 1 & a_d  
\end{array}
\right),
\qquad 
A^*:=\mbox{diag}(\theta^*_0, \theta^*_1, \ldots, \theta^*_d).
\end{eqnarray*} 
We show the pair
$A,A^*$ is a Leonard pair on $V$.
To do this we apply
Definition
\ref{def:lprecall}.
Observe that $A$ is irreducible tridiagonal and $A^*$
is diagonal. Therefore 
condition (i) of 
Definition
\ref{def:lprecall} is satisfied by the basis for
$V$ consisting of the columns of $I$, where $I$ denotes
the identity matrix in
$\mbox{Mat}_{d+1}(\K)$.
To verify condition (ii) of
Definition
\ref{def:lprecall},
we display an invertible matrix
$X$
such that $X^{-1}AX$ is diagonal and 
$X^{-1}A^*X$ is  irreducible tridiagonal.
Let 
$X$ denote the matrix in
$\mbox{Mat}_{d+1}(\K)$ 
that has entries
\begin{eqnarray}
X_{ij}&=&\frac{p_i(\theta_j)p^*_j(\theta^*_0)}
{x_1x_2\cdots x_i}
\label{eq:xver1}
\\
&=&\frac{p^*_j(\theta^*_i)p_i(\theta_0)}
{x_1x_2\cdots x_i}
\label{eq:xver2}
\end{eqnarray}
$0 \leq i ,j\leq d$. 
The matrix $X$ is invertible since it is essentially Vandermonde.
Using 
(\ref{eq:rec2introS99})
and 
(\ref{eq:xver1})
we find
$AX = XH$ where
$H=\mbox{diag}(\theta_0, \theta_1, \ldots, \theta_d)$.
Apparently $X^{-1}AX$ is equal to $H$ and is therefore diagonal.
Using (\ref{eq:rec3introS99})
 and
(\ref{eq:xver2})
we find
$A^*X=XH^*$ where
\begin{eqnarray*}
\label{eq:asmonmat}
H^*:=\left(
\begin{array}{ c c c c c c}
a^*_0 & x^*_1  &      &      &   &{\bf 0} \\
1 & a^*_1  &  x^*_2   &      &   &  \\
  & 1  &  \cdot    & \cdot  &   & \\
  &   & \cdot     & \cdot  & \cdot   & \\
  &   &           &  \cdot & \cdot & x^*_d \\
{\bf 0} &   &   &   & 1 & a^*_d  
\end{array}
\right).
\end{eqnarray*} 
Apparently $X^{-1}A^*X$ is equal to $H^*$ and is therefore irreducible
tridiagonal.
Now condition (ii) of Definition
\ref{def:lprecall} is satisfied by the basis for $V$ consisting
of the columns of $X$.
 We have now shown the pair $A,A^*$
 is a Leonard pair  on $V$.
Pick an integer $j$ $(0 \leq j \leq d)$.
Using $X^{-1}AX = H$ we find
 $\theta_j$ is the eigenvalue of 
$A$ associated with column $j$ of $X$.
From the definition of $A^*$ we find
 $\theta^*_j$ is the eigenvalue of 
$A^*$ associated with column $j$ of $I$.
Let $E_j$ (resp. $E^*_j$)
denote the primitive idempotent of $A$ (resp. $A^*$)
for $\theta_j$ (resp. $\theta^*_j$).
From our above comments the sequence
$\Phi:=(A; A^*; \lbrace E_i\rbrace_{i=0}^d; \lbrace E^*_i\rbrace_{i=0}^d)$
is a Leonard system. From the construction $\Phi$ is over $\K$.
We show
(\ref{eq:172L})
is the $MPS$
of $\Phi$.
To do this is suffices to show
 $a_i=\mbox{tr}(E^*_iA)$ for 
$0 \leq i \leq d$ and 
 $x_i=\mbox{tr}(E^*_iAE^*_{i-1}A)$ for $1 \leq i \leq d$.
Applying 
Lemma \ref{lem:bmat}(i),(ii) to $\Phi$
(with $v_i=\mbox{column  $i$ of   $I$}$, $B=A$)
we find
 $a_i=\mbox{tr}(E^*_iA)$ for 
$0 \leq i \leq d$ and 
 $x_i=\mbox{tr}(E^*_iAE^*_{i-1}A)$ for $1 \leq i \leq d$.
Therefore
(\ref{eq:172L})
is the $MPS$
of $\Phi$.
We show
(\ref{eq:172R})
is the dual $MPS$
of $\Phi$.
Applying
Lemma \ref{lem:bmat}(i),(ii) to $\Phi^*$
(with $v_i=\mbox{column  $i$ of $X$}$, $B=H^*$)
we find
 $a^*_i=\mbox{tr}(E_iA^*)$ for 
$0 \leq i \leq d$ and 
 $x^*_i=\mbox{tr}(E_iA^*E_{i-1}A^*)$ for $1 \leq i \leq d$.
Therefore 
(\ref{eq:172R})
is the dual $MPS$
of $\Phi$.
From the construction
we find (\ref{eq:173L}) (resp.  
(\ref{eq:173R}))
is the eigenvalue sequence
(resp.
dual eigenvalue sequence) of $\Phi$.
We show $\Phi$ is uniquely determined
by
(\ref{eq:172L})--(\ref{eq:173R})
up to isomorphism of Leonard systems.
Recall that $\Phi$ is determined
up to isomorphism of Leonard systems
by its own parameter array. We show the parameter array of 
$\Phi$ is determined by
(\ref{eq:172L})--(\ref{eq:173R}).
Recall the parameter array consists of
the eigenvalue sequence, the dual eigenvalue sequence,
the first split sequence and the second split sequence.
We mentioned earlier that the eigenvalue sequence of
$\Phi$ is
(\ref{eq:173L}) and the dual eigenvalue sequence
of $\Phi$ is 
(\ref{eq:173R}).
By Lemma \ref{lem:p0} the first split sequence of
$\Phi$ is determined by
(\ref{eq:172L})--(\ref{eq:173R}).
By this and
Theorem \ref{thm:xform} we find the second split sequence of 
$\Phi$ is determined by
(\ref{eq:172L})--(\ref{eq:173R}).
We have now shown the parameter array of $\Phi$ is
determined by 
(\ref{eq:172L})--(\ref{eq:173R}). We now see
that $\Phi$ 
is uniquely determined by
(\ref{eq:172L})--(\ref{eq:173R})
up to isomorphism of Leonard systems.
\hfill $\Box $ \\

\section{Acknowledgment} 
The author would like to thank Brian Hartwig and Darren Neubauer 
for giving this paper a close reading and offering many valuable
suggestions.

\noindent Paul Terwilliger \hfil\break
Department of Mathematics \hfil\break
University of Wisconsin \hfil\break
480 Lincoln Drive \hfil\break
Madison, Wisconsin, 53706 USA 
\hfil\break
email: terwilli@math.wisc.edu \hfil\break

\end{document}